\documentclass[review,authoryear]{elsarticle}

\usepackage{hyperref}
\usepackage{amsmath}
\usepackage{breqn}
\usepackage{calc}  
\usepackage{enumitem}  
\usepackage{amsfonts}
\usepackage{mathrsfs}
\usepackage{algpseudocode}
\usepackage{xcolor}
\usepackage{algorithm}
\usepackage{bigstrut}
\usepackage{rotating}
\usepackage[top=3.5cm,bottom=3.5cm,left=3.5cm,right=3.5cm]{geometry}
\usepackage{subcaption}
\usepackage{graphicx, amsfonts,amsmath , boldline, epstopdf,amssymb }
\usepackage{tabu, multirow,multicol ,booktabs ,setspace ,algcompatible ,mathrsfs,dsfont, array, boldline, makecell, booktabs }
\usepackage{url}
\usepackage{epstopdf}
\usepackage{adjustbox, stackengine, color, colortbl, epsfig, cases, enumitem, mathtools, tikz, verbatim}
\usetikzlibrary{shapes,arrows.meta,decorations}
\definecolor{Abi}{rgb}{0.309803, 0.58039, 0.80392}
\definecolor{orange}{rgb}{1, 0.5019, 0}
\definecolor{Red}{rgb}{1, 0, 0}
\definecolor{orange}{rgb}{1, 0.50196, 0}
\definecolor{greenJ}{rgb}{0, 0.6590, 0.42}
\definecolor{Brown}{rgb}{0.588, 0.294, 0}
\definecolor{um}{rgb}{0.0824, 0.1294, 0.4196}
\definecolor{abikam}{rgb}{0.51, 0.93,  0.992}
\newcommand*{\Rectangle}[1][]{%
  \tikz[baseline=-0.6ex]\node[draw,minimum width=0.65cm,minimum height=2ex, fill={#1}]{};
}

\usepackage{hyperref}
\usepackage{amsmath}
\usepackage{breqn}
\usepackage{calc}  
\usepackage{enumitem}  
\usepackage{amsfonts}
\usepackage{mathrsfs}
\usepackage{algorithm}
\usepackage{multirow}
\usepackage{bigstrut}
\usepackage{booktabs}
\usepackage{rotating}
\usepackage{hyperref}
\usepackage{bm}
\usepackage[top=3.5cm,bottom=3.5cm,left=3.5cm,right=3.5cm]{geometry}
\usepackage{caption}
 
\usepackage{mathtools}
\usepackage[algo2e]{algorithm2e}

\usepackage{color}
\newcommand{\rev}[1]{{\color{black}#1}}

\DeclareMathOperator*{\argmin}{arg\,min}

\newtheorem{lemma}{Lemma}
\newtheorem{theorem}{Theorem}

\newtheorem{defn}{Definition}

\newcommand{\reals}{{I\kern-.35em R}}
\newcommand{\nats}{{I\kern-.35em N}}
\newcommand{\upto}{{\raise 1pt \hbox{$\scriptstyle \,\nearrow\,$}}}
\newcommand{\downto}{{\raise 1pt \hbox{$\scriptstyle \,\searrow\,$}}}
\newcommand{\eop}
	{\hfill{$\vcenter{\hrule height1pt \hbox{\vrule width1pt height5pt 
   	 \kern5pt \vrule width1pt} \hrule height1pt}$} \medskip}
\def\state #1. { \noindent{\bf#1.\enspace}}
\newcommand{\nargmaxinf}{\mathop{\rm argmaxinf}\nolimits}
\newcommand{\nargmin}{\mathop{\rm argmin}\nolimits}
\newcommand{\nargmax}{\mathop{\rm argmax}\nolimits}
\newcommand{\nsupinf}{\mathop{\rm supinf}\nolimits}

\newcommand{\zg}[1]{\textcolor{red}{${\textrm{Zhaomiao: }}${#1}}}

 \usepackage{caption}
 
\journal{Transportation Research: Part C}

\begin{document}

\begin{frontmatter}

\title{Intermediate Service Facility Planning in a Stochastic and Competitive Market: Incorporating Agent-infrastructure Interactions over Networks}



\author[mythirdaddress]{Sina Baghali}
\author[mythirdaddress]{Zhaomiao Guo\corref{mycorrespondingauthor}}
\author[mysecondaryaddress]{Julio Deride}
\author[mymainaddress]{Yueyue Fan}
\cortext[mycorrespondingauthor]{(Corresponding Author) Assistant Professor, Department of Civil, Environmental and Construction Engineering, Resilient, Intelligent, and Sustainable Energy Systems Cluster, University of Central Florida, Orlando, FL 32816. Phone: 407-823-6215, Email: guo@ucf.edu}

\address[mythirdaddress]{Department of Civil, Environmental and Construction Engineering\\
	University of Central Florida}

\address[mysecondaryaddress]{Department of Mathematics\\
	Universidad T\'ecnica Federico Santa Mar\'ia\\
	Santiago, Chile}

\address[mymainaddress]{Department of Civil and Environmental Engineering\\
	University of California
	Davis
}

\begin{abstract}
This paper presents a network-based multi-agent optimization model for the strategic planning of service facilities in a stochastic and competitive market.  We focus on the type of service facilities that are of intermediate nature, i.e., users may need to deviate from the shortest path to receive/provide services in between the users' planned origins and destinations. This problem has many applications in emerging transportation mobility, including dynamic ride-sharing hub design and competitive \rev{facility location and allocation} problems for alternative fuel vehicle refueling stations.  \rev{ The main contribution of this paper is establishing a new multi-agent optimization framework considering decentralized decision makings of facility investors and users over a transportation network and providing rigorous analyses of its mathematical properties, such as uniqueness and existence of system equilibrium. In addition, we develop an exact convex reformulation of the original multi-agent optimization problems to overcome computational challenges brought by non-convexity.} \rev{Extensive analysis on case studies showed how the proposed model can capture the complex interaction between different stakeholders in an uncertain environment. Additionally, our model allowed quantifying the value of stochastic modeling and information availability by exploring stochastic metrics, including value of stochastic solution (VSS) and expected value of perfect information (EVPI), in a multi-agent framework.}
   
\end{abstract}

\begin{keyword}
Intermediate Service Facility \sep Competitive Facility Location \sep Multi-agent Optimization \sep Convex Reformulation
\end{keyword}
\end{frontmatter}



\section{Introduction}


Facility location-allocation problems (FLAPs), which seek the best strategy for locating facilities and allocating demands to the facilities, have wide applications in transportation science, supply chain and logistics, and infrastructure systems \citep{cornuejols1983uncapicitated,melo2009facility, hekmatfar2009facility}. In this paper, we focus on the type of service facilities that are of intermediate nature, i.e., users may need to deviate from a predefined shortest path to receive/provide services in between the users' planned origins and destinations. \rev{The stakeholders we model include facility investors and facility users. Facility investors decide investment capacity and provide services to facility users to maximize their own profit. Facility users will make facility selection and routing decision to receive services.} This problem has many applications in emerging transportation mobility, including competitive facility location problems for alternative fuel vehicle refueling stations and dynamic ride-sharing hubs design. 





Despite variations tailored to specific domain applications, there are some common features shared in the facility planning for alternative fuel vehicles and emerging mobility that have not been systematically studied in a unified framework. First, service demand could appear at travel origins, destinations, and/or intermediate locations. For example, plug-in electric vehicles (PEVs) can charge at home, workplaces, or public charging stations; ride-sharing/crowdsourced drivers may prefer to pick up and drop off riders/goods with minimum deviation from their planned route. These behaviors require additional modeling capabilities to provide flexibility in capturing node-based or/and link-based demand with possible deviation from pre-defined paths in an endogenous manner.  The second feature is that the system involves multiple competitive stakeholders from both supply and demand sides, who are driven by self-interests. For example, individual facility providers invest in service facilities to maximize their own profits. Individual travelers accessing facilities aim to optimize their utility, which could include travel time and/or service costs/revenue. This feature requires a modeling framework that could capture different decision entities' interests and enable analysis at the system level, where performance is shaped collectively by all. The third feature concerns complex agent-infrastructure interactions, i.e., the close coupling between the users' choices of facilities and travel routes, facility providers' facility location decisions, and the resulting link travel time and locational service prices, that need to be studied over a network structure. 


This paper aims to develop a generalized modeling framework and efficient computational algorithms to model and analyze intermediate service facility planning (ISFP) in a competitive and stochastic market. Specifically, we make contributions in the following three aspects: (1) we propose a unified system modeling framework to model the decentralized decision-making on the supply and demand sides of a competitive and stochastic market, with rigorous analysis of its mathematical properties (including equilibrium existence and uniqueness); (2) we extended the Combined Distribution and Assignment (CDA) model \citep{Evans1976} to capture the coupling between route choice and intermediate facility choices over a congested transportation network; and (3) from a computational perspective, we tailored an exact convex reformulation to our proposed equilibrium model, which significantly improves the computation efficiency with guaranteed global convergence. 




The remainder of the paper is organized as follows. Section 2 discusses relevant literature regarding competitive ISFP. Section 3 presents the proposed modeling framework and the mathematical formulations. The solution properties and computational strategies are provided in Section 4. Numerical experiments on the Sioux Falls test network are presented in Section 5 to provide analytical and numerical insights.  Section 6 concludes the paper with a discussion and potential future extensions.  


\section{Literature Review}

Given the ever-growing body of literature on FLAPs, we will focus on the discussion of the literature from a perspective that highlights what distinguishes our work from the past studies.  The readers may refer to \citep{Owen1998, Hale2003, Snyder2006, Melo2009, daskin2011network} for more comprehensive reviews on classic facility location models.  Most facility location-allocation models are built with a central planner's perspective, assuming that the location choices and sizes of different facilities can all be controlled by a single spatial monopoly/planner. Covering \citep{farahani2012covering}, p-center \citep{lin2018p}, p-median \citep{hansen1997variable}, and flow-capturing \citep{hodgson1990flow} are classic in this category. In reality, however, an infrastructure system often involves multiple facility developers driven by self interests \citep{Plastria_01}.  

To capture competitive nature of the supply side, competitive FLAPs, pioneered by \cite{Hotelling_29}, have been proposed and developed in facility location literature\citep{hakimi_83, Eiselt_et_al_93, Miller_et_al_96, Aboolian_et_al_07, Friesz_07, Drezner_09, smith_09, kress_pesch_12}. \rev{A competitive FLAP concerns the problem of deciding the locations and/or capacity of competing facilities, such as shopping centers, charging stations, ride-sharing hubs, restaurants, and others. In contrast to the classic facility location problem, the configurations of competing facilities are decided by a set of competitors who aim to optimize their own benefits.} Competitive facility location models can be broadly categorized based on how competition (e.g. static/dynamic), demand (e.g. fix/elastic, discrete/continuous, deterministic/probabilistic), and decision space (e.g. discrete/network/continuous) are formulated. These studies typically focus on decentralized decision-making from the supply side while simplifying demand-side modeling. For example, most existing studies on competitive FLAPs consider nodal demand (i.e., demand appearing at discrete locations) \citep{Klose_Drexl_05} or demand continuous in space \citep{li2010continuum} without considering endogenous demand that could be influenced by the facility locations. In addition, they typically do not model the travel and routing behavior of facility users over transportation network. A few competitive FLAPs studies consider flow-based demand and user travel routes (e.g., \citep{berman1998flow, wu2003solving}), but adopt a central planner’s perspective, which may undermine the capability to forecast and analyze the facility network collectively shaped by multiple investors. In addition, even though studies \citep{Yang_Wong_00,Ouyang_et_al_15} have shown that transportation congestion and user's facility choice are closely coupled, existing studies typically consider congestion at the facility level \citep{guo2016infrastructure,luo2015placement} and assume exogenously given traffic congestion, travel routes, and facility service prices.

The flow-based FLAPs have been actively studied in the last decade, especially in the context of charging stations for EVs. Based on Flow Intercepting Location Model (FILM)  \citep{hodgson1990flow, berman1992optimal}, \cite{shukla2011optimization, wen2014locating} developed mathematical programming models to determine the cost-effective charging station locations to maximize the intercepted traffic flow. In those studies, potential travel path deviations are not considered. Several versions of FILM with detours have been proposed by \cite{berman1995locating}, including maximizing O-D flows intercepted subject to maximum detour allowance and minimizing total detours subject to covering all O-D flows. Building upon \cite{berman1995locating}, deviated paths were considered in \citep{ li2014heuristic, zockaie2016solving}. This school of literature may not explicitly consider the EV driving range. \cite{kuby2005flow} has proposed the flow refueling location model (FRLM) to take into account driving range limitations for alternative fuel vehicles, which are further developed in \citep{kuby2009optimization, lim2010heuristic, capar2012efficient, mirhassani2013flexible, kim2013network, de2017incorporating, wang2018siting, guo2018battery, he2018optimal, boujelben2019efficient}.

\rev{Intermediate FLAPs have been applied in the context of EVs en-route charging to faciliate EV adoption \citep{kchaou2021charging}. \cite{wang2019designing} designed a charging station capacity and location problem for intra-city travels of EVs and \cite{xu2020mitigate} developed a facility location problem for battery swapping of EVs to reduce the range anxiety of drivers. Both studies focused on fulfilling the charging requirements of drivers during their trips and did not consider the transportation network congestion and drivers' facility location and routing choice modeling. \cite{xu2020optimal2} considered the elastic demand of drivers along with their path deviation in FLAP. \cite{li2022optimal} proposed a metanetwork-based approach to model en-route charging station planing to improve the network-based algorithm in the branch-and-bound framework. A bi-level optimization approach is proposed by \cite{tran2021user} to locate charging stations by minimizing the total travel time and installation costs at the upper level and captures re-routing behaviours of travellers with their driving ranges at the lower level. Authors use an iterative algorithm to solve the bi-level problem which does not guarantee finding the global optimal point. \cite{schoenberg2022siting} developed a charging station siting and sizing problem with coordinated charging to facilitate both en-route and destination charging where  the focus is on the charging scheduling instead of transportation network modeling.} \rev{All of the above studies take a central planner’s perspective, where all charging facilities are deployed by a single decision-maker.}

Modeling decentralized decision-makers in facility location problems has gained more attention in recent studies \citep{guo2016infrastructure, zhao2020deployment,bao2021optimal,chen2020optimal}. For example, \cite{zhao2020deployment} studied the optimal location of new charging stations among the existing competitive stations to maximize the profit of private investors. In that study, probabilistic modeling is developed to model the decision-making of the drivers where the congestion of the transportation network did not play a role in  the charging station selection of the drivers. \cite{bao2021optimal} developed a bi-level problem for the optimal charging station location of en-route charging in congested networks. The authors assumed that charging prices are similar in all of the charging stations across all stations and do not influence users’ facility choice. \cite{chen2020optimal} proposed a similar bi-level optimization framework where an investor decides the facility locations and their capacity at the upper level, and drivers' choices are modeled at the lower level. Neither study considered path deviation and en-route charging, which are essential in this context.    

\cite{guo2016infrastructure} proposed a network-based multi-agent optimization modeling framework to explicitly capture the decentralized behaviors of multiple facility investors and users in the context of public fast charger planning. That study modeled service demand only at travel destinations in deterministic market conditions. In addition, the formulation was nonconvex in \citep{guo2016infrastructure}. This paper aims to generalize \citep{guo2016infrastructure} in the following three aspects. First, we relax the assumption that travelers receive facility services only at trip destinations by modeling both node- and flow-base facility service demand. Second, we consider a stochastic market where parameters, such as OD travel demand, link travel time, and operational costs, could be uncertain.These extended modeling capabilities improve the realism of the studied problem setting. Third, when the multi-agent optimization problem is coupled with high-dimensional stochastic parameters, the combined problem becomes too complex to be solved by solution approaches based on lopsided convergence of bivariate functions as proposed in \citep{guo2016infrastructure}. To overcome that challenge, we establish an exact convex reformulation of the proposed modeling framework, which leads to significantly improved computational efficiency.

\section{Methodology}\label{sec:meth}

\subsection{Problem Description and Modeling Framework}

Our goal is to investigate the long-term equilibrium patterns of intermediate service facilities, considering the interactions between stakeholders from both facility supply and demand sides. On the facility supply side, we consider multiple investors, each of whom makes facility deployment and operational decisions to maximize its own profits. We assume each facility provider does not have the sufficient market power to strategically influence the locational service prices through its own decision-making (i.e., service providers are perfectly competitive)\footnote{We acknowledge that some markets may not fall into the perfect competition category, such as US electricity wholesale market, where the entry barriers and capacity constraint may lead to imperfect competition, especially during contingency \citep{guo2017stochastic}. For those markets involving noticeable market power, an oligopolistic model, such as Cournot, Bertrand, or Hotelling model, would be more appropriate. These market settings are beyond the scope of this paper, and we shall leave the investigation of alternative market structures in the future.}. On the demand side, there are (many) potential facility service users who make individual choices both for facilities and travel routes in order to maximize their utilities, which may depend on facility service prices, locational preference, and the travel time to access the facility. Locational facility service prices and travel time are endogenously determined through the interactions between service supply and demand over the transportation network.

We model this problem in the framework of network-based multi-agent optimization problem with equilibrium constraints (N-MOPEC) \citep{guo2016infrastructure}, which reflects the ``selfish'' nature of each decision entity while simultaneously capturing the interactions among all over a complex network structure. MOPEC is originally proposed by \citep{ferris2013mopec}, which includes a wide variety of variational problems as special cases: variational inequalities, complementarity problems, fixed points problems, etc. MOPEC has wide applications in  economics \citep{Deride_et_al_15}, coupled transportation/power systems \citep{guo2021stochastic, baghali2022electric}, and ride-sourcing mobility systems \citep{afifah2022spatial}.

Consider a collection of agents $A$ whose decisions are denoted as $\boldsymbol{x}_A = (\boldsymbol{x}_a, a \in A)$.  A MOPEC model, in its general form, can be expressed as: 

\begin{equation}\label{eq:mopec_mop}
\boldsymbol{x}_a \in \textrm{argmax}_{\boldsymbol{x} \in  X_{\boldsymbol{p},\boldsymbol{x}_{-a}} \subset \reals^{n_a}}  \;  f_a(\boldsymbol{p}, \boldsymbol{x}, \boldsymbol{x}_{-a}), \; a \in A, 
\end{equation}
\rev{where $\boldsymbol{x}$ represents the vector of investor $a$'s decision variables and $X_{\boldsymbol{p},\boldsymbol{x}_{-a}}$ is the feasible set for the investor's problem which may depend on the system parameters $\boldsymbol{p}$ and other investors' decisions  $\boldsymbol{x}_{-a}$ ($-a$ means $A \setminus a$). $\reals^{n_a}$ represents the domain of the set with $n_a$ being the dimension of decision variables for investor $a$.  $f_a$ is agent $a$'s objective function depends on the decisions of the other agents  and system parameters ($\boldsymbol{p}$), which may be endogenously determined by the system, such as prices.} Parameters $\boldsymbol{p}$ and the decisions $\boldsymbol{x}_A$ resulting from the multi-optimization problem typically need to satisfy global equilibrium constraints, which can be formulated as a functional variational inequality (\ref{eq:mopec_ec}):

\begin{equation}\label{eq:mopec_ec}
D(\boldsymbol{p},\boldsymbol{x}_A) \in \partial g(\boldsymbol{p}), 
\end{equation}
where $g: \reals^d \rightarrow \overline{\reals}$ is a proper, lower semicontinuous and convex function and $D$ is a set-valued mapping from $\reals^d \times \reals^{\sum_{a\in A} n_a}$ to $\reals^d$. 

\begin{figure}[htbp]
\begin{center}
    \includegraphics[width=0.7\textwidth]{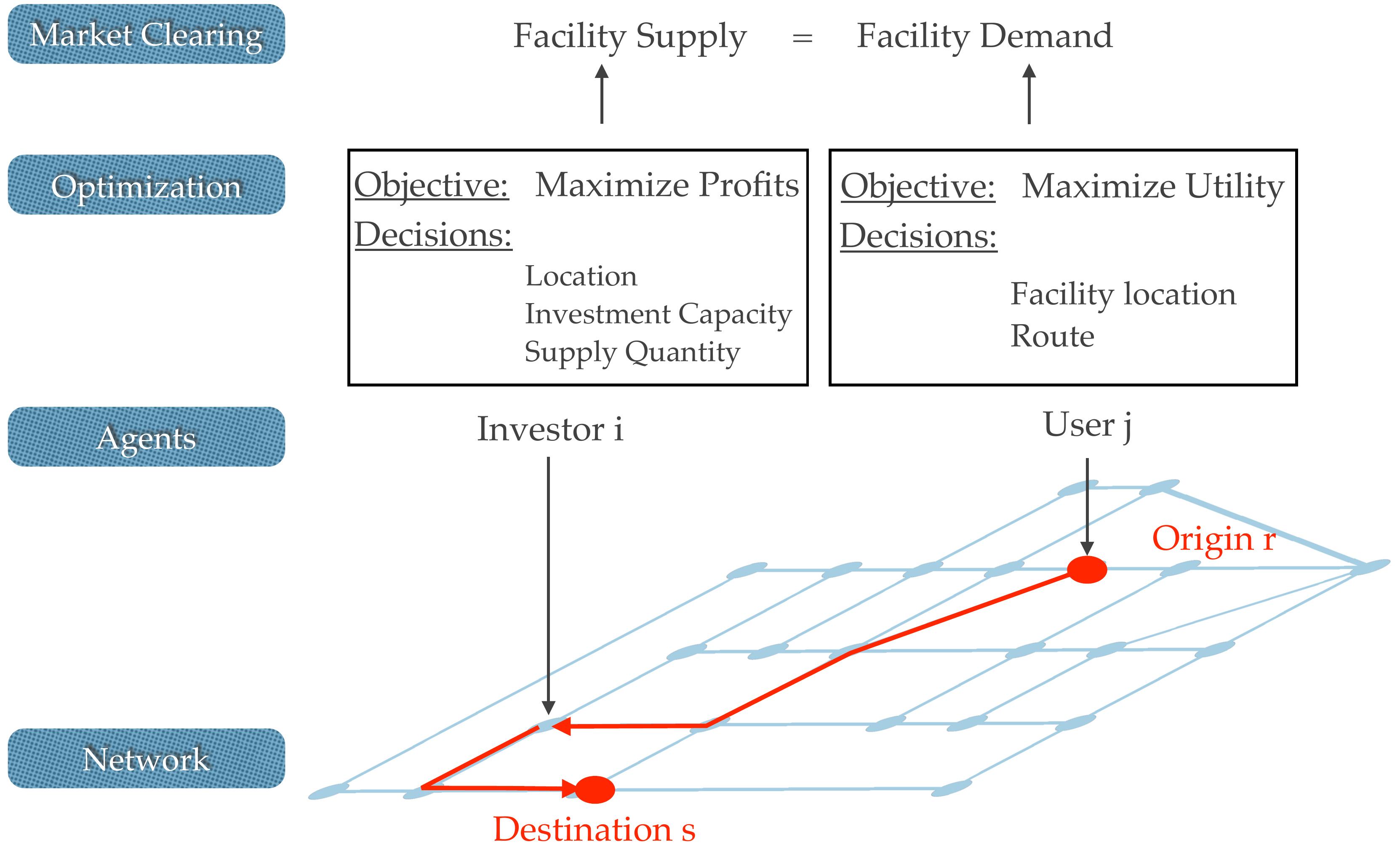}
\caption{Illustration of Network-based MOPEC}
\label{fig:mode_fram}
\end{center}
\end{figure}

The proposed N-MOPEC modeling framework in the context of competitive ISFP is illustrated in Figure \ref{fig:mode_fram}. We consider two categories of stakeholders: (1) individual investor $i$ ($\in \mathcal{I}$) decides the location, facility service capacity, and supply quantity to maximize his/her own profits; (2) individual service user $j$ ($\in \mathcal{J}$ ) travels from a specific origin and destination. User $j$ chooses facility service locations and travel routes to maximize his/her own utility.  Even though the decisions of these agents are made individually, they are interdependent due to the shared market, infrastructure, and resources. To ensure an equilibrium state is reached, market clearing conditions, i.e., supply equals demand at every facility location, also need to be imposed.  

\subsection{Detailed Formulation for Each Agent} 

\subsubsection{Modeling the Decisions of Facility Investors}\label{ssec:Inv}

Although this paper focuses on facility layout in the long run from investors' perspective, the effectiveness of planning decisions can not be properly evaluated without considering the performance in the operational stage. Therefore, we adopt a two-stage stochastic programming framework to distinguish between two types of decisions a facility investor has to make: (1) during the planning stage, each investor decides the capacities of facilities to invest facing future uncertainties, such as demand, access time, and marginal operational costs. (2) during the operational stage, uncertain parameters are revealed, and each investor will choose its supply quantities based on market locational prices and operational costs. Because of the assumption of a perfectly competitive market, without loss of generality, we can aggregate the decision-making of all investors into a representative one and use aggregate investment and operational cost functions to capture their collective decisions. The detailed formulation for the representative investor is presented in the model (\ref{eq:investor}). 


\begin{subequations}
\begin{align}
&
{\underset{c^k,g^{k}_{\boldsymbol{\xi}}\in\reals_{+},k\in K,\boldsymbol{\xi}\in\boldsymbol{\Xi}}{\text{maximize}}}
&&\hspace{-1cm}\mathbb{E}_{\boldsymbol{\xi}}\sum_{k \in K}\left[ \rho^{k}_{\boldsymbol{\xi}}g^{k}_{\boldsymbol{\xi}} - \phi_g( g^{k}_{\boldsymbol{\xi}})\right]-\sum_{k\in K}\phi_c(c^k)\label{obj:inve}\\
&\makebox[3cm]{\text{subject to}}
&&\hspace{-1cm}g_{\boldsymbol{\xi}}^{k} - c^k \leq 0,\forall k \in K,\boldsymbol{\xi}\in\boldsymbol{\Xi}\label{cons:capa}
\end{align}
\label{eq:investor}
\end{subequations}
where: 
\begin{description}[leftmargin=!,labelwidth=\widthof{12345}]
\itemsep-0.7em 
\item[$K$] set of candidate investment locations, indexed by $k$; 
\item[$\boldsymbol{\Xi}$] vector set of uncertain parameters, indexed by $\boldsymbol{\xi}$;
\item[$c^k$] investment capacity allocated at location $k$; 
\item[$g_{\boldsymbol{\xi}}^{k}$] total supply at location $k$ in scenario $\boldsymbol{\xi}$; 
\item[$\rho_{\boldsymbol{\xi}}^k$] unit service price at location $k$ in scenario $\boldsymbol{\xi}$, endogeneously determined by the market; 
\item[$\mathbb{E}_{\boldsymbol{\xi}}$] expectation with respect to the uncertain parameters $\boldsymbol{\xi}$; 
\item[$\phi_c (\cdot)$] aggregate capital cost function with respect to facility capacity; 
\item[$\phi_g (\cdot)$] aggregate operational cost function with respect to supply quantity.
\end{description}


The objective function (\ref{obj:inve}) maximizes the expected net profits, calculated as the expectation of the total revenues $\sum_{k \in K} \rho^{k}_{\boldsymbol{\xi}}g^{k}_{\boldsymbol{\xi}}$  minus the operating cost $\sum_{k \in K} \phi_g  ( g^{k}_{\boldsymbol{\xi}})$, minus the total investment cost during planning stage $\sum_ {k \in K} \phi_c (c^k )$. Total investment costs could include costs associated with land acquisition, construction, and equipment purchase. Constraint (\ref{cons:capa}) is the capacity constraint that ensures the supplied quantities at each location $k$ and scenario $\xi$ do not exceed its total capacity. The remaining constraints are non-negative restrictions. Note that throughout the paper, we denote vectors in lowercase bold font.


The interpretation of uncertainties $\boldsymbol{\xi}$ is two-fold. First, investors can not predict the future service demand due to uncertain factors, such as total demand (e.g., EV adoption), market competition, travel/charging time/costs, etc. In this case, the interpretation of the probability of $\boldsymbol{\xi}$ is the probability of uncertain parameters. Second, the state of the systems (e.g., facility service demand) may change over time, which can be grouped into homogeneous time segments (e.g., peak and off-peak hours). In this case, the probability of $\boldsymbol{\xi}$ measures the duration percentage of certain homogeneous time segments in the studied horizon.  In other words, $\boldsymbol{\xi}$ can represent a realization of uncertain parameters and/or a specific homogeneous time segment. 


\rev{Note that $\phi_c(\cdot)$ and $\phi_g(\cdot)$ are aggregate capital and operational cost functions at each location.} In this paper, we assume $\phi_c^k (\cdot)$ and $\phi_g^k (\cdot)$ to be convex functions, e.g., linear function or a quadratic form with positive leading coefficients.  Besides mathematical convenience, a convex production cost function implies two desired properties: (1) as service demand at a location increases, it may cause congestion in the upstream supply chain, which leads to a higher marginal cost; (2) as demand increases, higher-cost production resources may start to be utilized, due to capacity limitation of lower cost resources. \rev{For example, due to space limitations, the earlier investment can be made at a location with cheaper rent and/or construction costs. However, later investment may have to be built at a more expensive location. Capacity cost functions with increasing marginal costs are widely used to model the cost  of charging station \cite{ghamami2020refueling, ghamami2016general, 	guo2018battery}.  In terms of operational costs, when the facility needs certain resources to operate and the resource supply are limited, the marginal production cost is usually monotone increasing with production quantity. For example, for charging facility, the energy prices increase with the demand quantity because cheaper energy resources will be dispatched first. For shared mobility, to attract additional unit of drivers, transportation network companies usually need to pay higher prices.}

\subsubsection{Modeling the Decisions of Facility Users in a Congested Transportation Network}

Facility users' behaviors (facility choice and route choice) are affected by not only the characteristics of facilities but also the transportation network.  The combined Distribution and Assignment (CDA) model (see, e.g., \citep{Sheffi_85, lam_huang_92}) has been demonstrated effective in terms of integrating discrete choices (e.g., mode choices and destination choices) and traffic assignment in the context of charging infrastructure planning \citep{he2013optimal, guo2016infrastructure}.  In this study, instead of restricting the service location to be at the travel destinations, we propose a Generalized Combined Distribution and Assignment (GCDA) model, in which the facility location can be at either origin, destination, or anywhere in between. Since all decisions from the demand side are operational decisions and scenario dependent, we omit the notation $\boldsymbol{\xi}$ for brevity throughout this subsection.  

A multinomial logit model is used to describe the choice of different facility locations $k$ from origins $r$ to destination $s$, with the utility function defined in (\ref{eq:utility}). 

\begin{equation} \label{eq:utility}
U^{rsk} = \beta_0^k -\beta_1 t^{rsk} - \beta_2 \rho^k e^{rs} + \epsilon^{rsk}
\end{equation}
where: 
\begin{description}[leftmargin=!,labelwidth=\widthof{12345}]
\itemsep-0.7em 

\item[$U^{rsk}$]: utility measure of a user to go from $r$ to $s$ and receive service at $k$;
\item[$\beta$] : utility function parameters (model input);
\item[$t^{rsk}$] : equilibrium travel time from $r$ to $s$, with detour to service location $k$; 
\item[$e^{rs}$] : average service demand from $r$ to $s$ (model input);
\item[$\epsilon^{rsk}$] : error term of utility from $r$ to $s$, with detour to service location $k$. $\epsilon^{rsk}$ follows extreme value distribution.
\end{description}

The utility function of a traveler from origin node $r$ to destination node $s$ choosing facility $k$ is assumed to be the summation of four parts: locational specific attractiveness factor ($\beta_0$), travel time ($\beta_1$), and service cost ($\beta_2$), and error term ($\epsilon^{rsk}$). \rev{$\beta_2 \rho^k e^{rs}$ in the utility function  (\ref{eq:utility}) is monetary service costs. $\rho^k e^{rs}$ represents the price that drivers pay for the service at facility $k$. This price is calculated as the unit cost of service at station $k$ (i.e. $\rho^k$) multiplied by the average charging quantity ($e^{rs}$). $\beta_2$ is a utility coefficient representing the disutility for each unit of money spent.} \rev{Here, we have considered vehicles to have similar facility demand $e^{rs}$ coming from each O-D pair. However, this assumption can be easily relaxed by categorizing drivers from each origin node based on their different levels of facility needs. In case of EVs, for example, the $e^{rs}$ can be categorized into different homogeneous groups to model EVs with different charging demand levels.}  Different exogenous utility factors can be included in the utility function without affecting the key modeling and computational strategies proposed in this study.  \rev{In addition, although service time is not explicitly modeled in the utility function (\ref{eq:utility}), in our model one may consider adding a dummy link connecting the closest transportation node to the facility. The travel time of this dummy link represents the service time needed at that facility depending on the service capacity.}


Denote the transportation network by a directed graph $\mathcal{G} = (\mathcal{N}, \mathcal{A})$ , where $\mathcal{N}$  is the set of nodes (indexed by $n$) and  $\mathcal{A}$ is the set of links (indexed by $a$).  A node can represent  a TAZ (source/sink of aggregated travel demand), a transport hub, or an intersection. A link can represent a path or a physical road section that connects two nodes. The GCDA model for a scenario $\xi (\in \boldsymbol{\Xi})$ is formulated in (\ref{eq:cda}).

\begin{subequations}
\begin{align}
& \underset{\boldsymbol{\hat{x}},\boldsymbol{\check{x}}, \boldsymbol{x},  \boldsymbol{q} \geq \boldsymbol{0}}{\text{minimize}}
& & & &\sum_{a \in \mathcal{A}} \int_{0}^{v_a} t_a(u) \mathrm{d}u \nonumber \\
& & & & &+ \frac{1}{\beta_1}\sum_{r \in R}\sum_{s \in S}\sum_{k \in K^{rs}} q^{rsk}\left(\ln q^{rsk} - 1 + \beta_2\rho^ke^{rs} - \beta_0^k\right) \label{obj:CDA}\\
& \text{subject to}
& & & &v_a = \sum_{r \in R} \sum_{s \in S} \sum_{k \in K^{rs}} (\hat{x}_{a}^{rsk} + \check{x}_{a}^{rsk}), \forall a \in \mathcal{A} \label{cons:v_x}\\
& & & (\boldsymbol{\gamma}) & &\boldsymbol{\hat{x}}^{rsk} + \boldsymbol{\check{x}}^{rsk} = \sum_{p \in P^{rsk}}(B_{\hat{p}} + B_{\check{p}})x_p , \; \forall r \in R, s \in S, k \in K^{rs}\label{cons:x_p}\\
& & & (\boldsymbol{\hat{\lambda}})& &  A\boldsymbol{\hat{x}}^{rsk} = q^{rsk} E^{rk}, \; \forall r \in R, s \in S, k \in K^{rs}\label{cons:x1_q}\\
& & & (\boldsymbol{\check{\lambda}}) & &  A\boldsymbol{\check{x}}^{rsk} = q^{rsk} E^{ks}, \; \forall r \in R, s \in S, k \in K^{rs}\label{cons:x2_q}\\
& & & (\mu^{rs}) & & \sum_{k \in K^{rs}} q^{rsk} = d^{rs}, \forall r \in R, s \in S\label{cons:q_d}
\end{align}
\label{eq:cda}
\end{subequations}
where: 
\begin{description}[leftmargin=!,labelwidth=\widthof{12345}]
\itemsep-0.7em
\item[$v_a$] : traffic flow on link $a$; 
\item[$t_a(\cdot)$] : travel time function of link $a$, e.g. the Bureau of Public Roads (BPR) function; 
\item[$d^{rs}$] : travel demand from $r$ to $s$ (model input); 
\item[$q^{rsk}$] : traffic flow from $r$ to $s$ and service at $k$;
\item[${x}_p$] : traffic flow on path $p$;
\item[$\hat{x}_a^{rsk}$] : traffic flow on link $a$ that belongs to the travel from $r$ to $k$ associated with Origin-Service-Destination triple $rks$. $\boldsymbol{\hat{x}}^{rsk}$ represents the vector form of $\hat{x}_a^{rsk}$ for all links;
\item[$\check{x}_a^{rsk}$] : traffic flow on link $a$ that belongs to the travel from $k$ to $s$ associated with Origin-Service-Destination triple $rks$. $\boldsymbol{\check{x}}^{rsk}$ represents the vector form of $\hat{x}_a^{rsk}$ for all links;
\item[$A$] : node-link incidence matrix of network, with $1$ at starting node and $-1$ at ending node;
\item[$\hat{p}$] : sub-path of path $p \in P^{rsk}$ that connect $r$ to $k$;
\item[$\check{p}$] : sub-path of path $p \in P^{rsk}$ that connect $k$ to $s$;
\item[$B_p$] : link-path $p$ incidence vector, with $i$th row equals to $1$ if path $p$ includes link $i$ and $0$ otherwise;
\item[$E^{ij}$] : O-D incidence vector of O-D pair $ij$ with $1$ at origin $i$, $-1$ at destination $j$;
\item[$\gamma, \check{\lambda}, \hat{\lambda}, \mu$] : dual variables of the corresponding constraints. 
\end{description}

Constraint (\ref{cons:v_x}) calculates the aggregate link flow $v_a$ from the link flow associated with $rsk$: $\hat{x}_a^{rsk}$ and $\check{x}_a^{rsk}$; Constraint (\ref{cons:x_p}) guarantees there is always a feasible path flow solution $x_p (p \in P)$ that can yield a given link flow pattern;  Constraint (\ref{cons:x1_q}, \ref{cons:x2_q}) ensures the flow conservation at each node, including the origin, intermediate stop, and destination nodes; Constraint (\ref{cons:q_d}) guarantees the sum of the flow to all facilities equals to the total travel demand between each origin-destination pair. Note that the OD demand that does not need access to facility service can be considered as background traffic in model (\ref{eq:cda}). In addition, the total demand $d^{rs}$ may be elastic, and our modeling framework can be naturally extended to consider elastic travel demand depending on travel distance, time, service congestion, and costs. Those who are interested in elastic demand can refer to \citep{Berman_Kaplan_87, Aboolian_et_al_12, Berman_Drezner_06}. The rest of the constraints set non-negative restrictions on path/link flow and trip distribution. 

In the objective function (\ref{obj:CDA}), the first term corresponds to the total user cost as modeled in a conventional static traffic equilibrium model; the second term involving $q\ln q$ corresponds to the entropy of trip distribution, and the remaining terms correspond to the utility measure (\ref{eq:utility}) of the travelers.  This objective function does not have a physical interpretation, but it guarantees the first Wardrop principle \citep{Wardrop_52} and the multinomial logit facility choice assumption being satisfied, as formally stated in {\bf Lemma \ref{lem:GCDA}}. 

\begin{lemma}{(Generalized Combined Distribution and Assignment)}\label{lem:GCDA}
The optimal solutions $(\boldsymbol{\hat{x}^\ast},\boldsymbol{\check{x}^{\ast}},$ $ \boldsymbol{x^\ast},  \boldsymbol{q^\ast} )$ of  problem (\ref{eq:cda}) are the equilibrium solutions for the service location choice with logit facility demand functions and Wardrop user equilibrium.
\end{lemma}  

\state Proof. See \ref{app:pfs}. 
\eop

\begin{figure}[htbp]
\begin{center}
    \includegraphics[width=0.7 \textwidth]{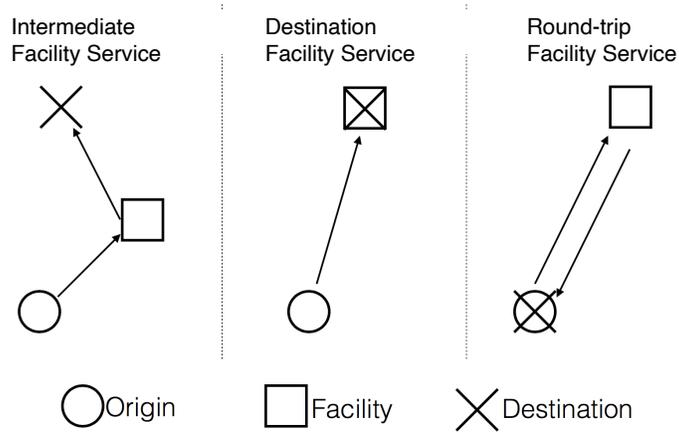}
\caption{Special Cases of GCDA Model}
\label{fig:special_case}
\end{center}
\end{figure}

The GCDA model proposed here can include three special cases, denoted as {\it Intermediate Facility Service }, {\it Origin/Destination Facility Service}, and {\it Round-trip Facility Service}, see Figure \ref{fig:special_case}. {\it Intermediate Facility Service} case represents when facility users access the facility service on their way to a destination, such as refueling, banking, and convenience store services.  In the {\it Destination Facility Service} case, drivers choose their destinations and, in the mean time, receive service at their destinations. For example, EV drivers may choose a restaurant and charge their vehicles at the same time. This case is precisely the conventional CDA model\citep{Sheffi_85}. In the {\it Round-trip Facility Service} case, a user who starts from the origin will make a dedicated trip to a facility location and will need to go back to the same origin after receiving service. For example, employees who have lunch and need to go back to the workplace afterward. To degenerate our model to {\it Destination Facility Service} case, we can specify $s = k$, while {\it Round-trip Facility Service} can can be incorporated in our model by specifying $r = s$.

\rev{Note that the framework is not a stochastic user equilibrium (SUE) per the definition by \cite{daganzo1977stochastic} and we do not need to generate a predefined set of paths in advanced. The ``stochasticity'' in this paper refers to the uncertainties facility investors faced when they make long-term planning decisions. Model (\ref{eq:cda}) is scenario dependent and we aims to model the equilibrium traffic patterns given each specific realization of scenario. In other words, travelers do not face uncertainties when they make facility and routing choices. Model (\ref{eq:cda}) is an extension of standard Beckmann's formulation of Wardrop user equilibrium \citep{Wardrop_52} by considering intermediate facility location choice using Logit model. A closely related category of model is called combined distribution and assignment (CDA) model \citep{Sheffi_85, lam_huang_92}, where travelers choose destinations instead of intermediate facility locations. }

\subsubsection{Market Clearing Conditions}\label{sec:mark_clea}

Lastly, to ensure the market is stabilized, we formulate the market clearing conditions that requires the total demand ($\sum_{r \in R}\sum_{s \in S}e^{rs} q^{rsk}_{\boldsymbol{\xi}}$) equal the total supply ($g_{\boldsymbol{\xi}}^{k}$) at each facility location in each scenario, as in (\ref{cons:mark_clea}). \rev{Recall that $\boldsymbol{q}^{rsk}$ represents the travel demand of $rs$ that receive services at location $k$. $\boldsymbol{e}^{rs}$ represents the average quantity of service demand per user from $rs$. Therefore, the total service demand is calculated as $\sum_{r \in R}\sum_{s \in S}e^{rs} q^{rsk}_{\boldsymbol{\xi}}$.} 


Locational service prices $\rho^k_{\boldsymbol{\xi}}$ can be interpreted as the dual variables of (\ref{cons:mark_clea}). Locational service prices  $\rho^k_{\boldsymbol{\xi}}$ influence the decision making of both supply and demand sides (i.e., problems (\ref{eq:investor}) and (\ref{eq:cda})) to optimize their own objectives. On the other hand, supply and demand imbalance will influence the locational service prices. We focus on estimating the locational service prices that can lead to a market clearing in an equilibrium state. Note that due to network congestion and accessibility cost, the prices may vary by location even if the services offered at each location are identical.

\begin{equation} 
(\rho^k_{\boldsymbol{\xi}}) \ \ \  g_{\boldsymbol{\xi}}^{k} - \sum_{r \in R}\sum_{s \in S}e^{rs} q^{rsk}_{\boldsymbol{\xi}}= 0, \; \forall k \in K, \boldsymbol{\xi} \in \boldsymbol{\Xi}.
\label{cons:mark_clea}
\end{equation}


\subsubsection{System Equilibrium }\label{ssec:systeq}
 
The decisions of all participants in this system are interdependent and should be modeled and solved simultaneously as a whole system. Following the notion of Nash equilibrium, at system equilibrium, a unilateral decision change of one agent given the market clearing price ${\rho}$ and other agents' decisions would diminish his/her pay-off. We state the system equilibrium more formally in Definition \ref{defn:Syst_Equi}.

\begin{defn}\label{defn:Syst_Equi}
{\it(system equilibrium)}. {\normalfont The equilibrium state of the system is that all facility providers achieve their own optimality (i.e., problem (\ref{eq:investor})) and facility users achieve their optimality  (i.e., problem (\ref{eq:cda})), given market clearing price $\boldsymbol{\rho}$ and all other agents' decisions.  In addition, the market at each location is cleared by condition (\ref{cons:mark_clea}).} Thus, a system equilibrium is defined by an investor strategy $(\boldsymbol{c}^*,\boldsymbol{g}^*)$, a GCDA traffic pattern $({\boldsymbol{\hat{x}}}^*,\boldsymbol{\check{x}}^*,\boldsymbol{x}^*,\boldsymbol{q}^*)$, and a vector of prices $\boldsymbol{\rho}^*$, such that
\begin{eqnarray}
(\boldsymbol{c}^*,\boldsymbol{g}^*)\,\mbox{solves the Investor problem for price }\boldsymbol{\rho}^*\mbox{ in }\, (\ref{eq:investor}),\\
(\boldsymbol{\hat{x}}^*_{\boldsymbol{\xi}},\boldsymbol{\check{x}}^*_{\boldsymbol{\xi}},\boldsymbol{x}^*_{\boldsymbol{\xi}},\boldsymbol{q}^*_{\boldsymbol{\xi}})\,\mbox{solves each the GCDA problem for price }\boldsymbol{\rho}^*_{\boldsymbol{\xi}}\mbox{ in }\,  (\ref{eq:cda}), \forall \boldsymbol{\xi} \in \boldsymbol{\Xi}\\
\mbox{and }(\boldsymbol{g}^*,\boldsymbol{q}^*)\mbox{ satisfy the market clearing conditions in } (\ref{cons:mark_clea}).
\end{eqnarray}
\end{defn}



\section{Solution Methods}\label{sec:solu_meth}

The proposed modeling framework in Section \ref{sec:meth} is a highly non-convex problem due to the complementarity conditions, which is challenging to solve especially with large number of scenarios. In this section, we propose an exact convex reformulation for the original N-MOPEC problem, which can lead to further scenario decomposition for scalability. In addition, leveraging convex reformulation, we prove the existence and uniqueness of the equilibrium solution to the original  N-MOPEC problem.

\subsection{Exact Convex Reformulation}\label{ssec:centralized}
Solving the system equilibrium state per Definition \ref{defn:Syst_Equi} is non-trivial due to the complementary nature of the model formulation. We proposed an exact convex reformulation to solve the system equilibrium, which is shown in model \eqref{eq:combined}.  A similar approach is described in \citep{Dv20StoESO} and has been  proposed to convexify multi-agent system equilibrium problems in coupled transportation and power systems \citep{guo2021stochastic, baghali2022electric} and ride-sourcing systems \citep{afifah2022spatial}. In the proposed exacted convex reformulation (i.e., model (\ref{eq:combined})), the objective function is to minimize the combined social non-transactional costs, i.e., the linear combination of the investors' costs and the normalized\footnote{CDA objective is normalized by $\beta_1/\beta_3$. The intuition of this step is to convert the unit of CDA objective to \$. More rigorous proof of why we adopt this form can be seen in Lemma \ref{lem:equivalency}.} CDA objective function without considering the terms associated with the price vector $\boldsymbol{\rho}$. The constraints of model (\ref{eq:combined}) include the constraints of investors~\eqref{cons:capa}, GCDA~\eqref{cons:v_x}-\eqref{cons:q_d}, and the market clearing condition~\eqref{cons:mark_clea}. 

\begin{subequations}
\begin{align}
&
{\underset{(\boldsymbol{c,{g},{\hat{x}},{\check{x}}, {x},  {q}})}{\text{minimize}}}
&&\hspace{-1cm}\sum_{k\in K}\phi_c(c^k)+ \mathbb{E}_{\boldsymbol{\xi}}\sum_{k \in K}\left[\phi_g( g^{k}_{\boldsymbol{\xi}})\right]\nonumber \\
& & & \hspace{-1cm}+ \mathbb{E}_{\boldsymbol{\xi}}\left\{\frac{\beta_1}{\beta_3}\sum_{a \in \mathcal{A}} \int_{0}^{{v}_{a,\boldsymbol{\xi}}} t_a(u) \mathrm{d}u + \frac{1}{\beta_3}\sum_{r \in R}\sum_{s \in S}\sum_{k \in K^{rs}} {q}^{rsk}_{\boldsymbol{\xi}}\left(\ln {q}^{rsk}_{\boldsymbol{\xi}} - 1 - \beta_0^k\right)\right\} \label{obj:centralized}\\
&\makebox[3cm]{\text{subject to}}
&(\lambda_{\boldsymbol{\xi}}^k)&\hspace{0.2cm}g_{\boldsymbol{\xi}}^{k} - \sum_{r \in R}\sum_{s \in S}e^{rs} q^{rsk}_{\boldsymbol{\xi}}= 0, \; \forall k \in K, \boldsymbol{\xi} \in \boldsymbol{\Xi}.\label{con:cceq}\\
& & & (\boldsymbol{c,{g}})\,\mbox{satisfies\, constraint~\eqref{cons:capa}}\\
& & & (\boldsymbol{\hat{x}_{\xi}, \check{x}_{\xi}, x_{\xi},  q_{\xi}})\,\mbox{satisfies\, constraints~\eqref{cons:v_x}-\eqref{cons:q_d}}_{\boldsymbol{\xi}},\,\forall \boldsymbol{\xi}\in \boldsymbol{\Xi}\label{con:cccda}
\end{align}
\label{eq:combined}
\end{subequations}


To investigate the existence and uniqueness of system equilibrium and its relationship with the solutions of model \eqref{eq:combined}, we first prove that model \eqref{eq:combined} is a strictly convex optimization problem under mild conditions (see Lemma \ref{lem:convex_existence}). Then we show that any solutions solving model \eqref{eq:combined} will also satisfy the system equilibrium definition \ref{defn:Syst_Equi} (see Lemma \ref{lem:equivalency}). 

\begin{lemma}{(convexity of model \eqref{eq:combined} and solution uniqueness)}\label{lem:convex_existence}
	Model \eqref{eq:combined} is convex if $\phi_c(\cdot)$ and $\phi_g(\cdot)$ are convex functions and $t_a(\cdot)$ is monotone increasing. Furthermore, if $\phi_c(\cdot)$ and $\phi_g(\cdot)$ are strictly convex functions and $t_a(\cdot)$ is strictly monotone increasing, model \eqref{eq:combined} is strictly convex and has a unique solution.
\end{lemma} 

\state Proof. See \ref{app:pfs}. 

\begin{lemma}{(solutions of model (\ref{eq:combined}) and their relationship with system equilibria)}\label{lem:equivalency}
	Assume $\phi_c(\cdot)$ and $\phi_g(\cdot)$ are convex functions and $t_a(\cdot)$ is monotone increasing. $(\boldsymbol{c}^*,\boldsymbol{g}^*,\boldsymbol{\hat{x}}^*,\boldsymbol{\check{x}}^*, \boldsymbol{x}^*,  \boldsymbol{q}^*; \boldsymbol{\lambda})$ is a primal-dual solution of the model~\eqref{eq:combined} if and only if  $(\boldsymbol{c}^*,\boldsymbol{g}^*,\boldsymbol{\hat{x}}^*,\boldsymbol{\check{x}}^*, \boldsymbol{x}^*,  \boldsymbol{q}^*)$ is a system equilibrium (Definition \ref{defn:Syst_Equi}), with equilibrium price vector $\rho_{\boldsymbol{\xi}}^k=\frac{\lambda_{\boldsymbol{\xi}}^k}{\pi_{\boldsymbol{\xi}}}$, for every $\boldsymbol{\xi}$ and $k$, where $\{\pi_{\boldsymbol{\xi}}:\boldsymbol{\xi} \in \boldsymbol{\Xi}\}$ is the probability distribution of $\boldsymbol{\xi}$.
\end{lemma} 

\state Proof. See \ref{app:pfs}. 

Lemma \ref{lem:equivalency} establishes the equivalency between our convex reformulation and the equilibrium model (i.e., model \eqref{eq:investor}, \eqref{eq:cda}, and \eqref{cons:mark_clea}). Furthermore, based on Lemma \ref{lem:convex_existence}, we can further discuss the uniqueness of the system equilibrium under strict convexity conditions of the exact reformulation, as shown in Theorem \ref{thm:exis_syst_equi}.

\begin{theorem}{(existence and uniqueness of system equilibrium)}\label{thm:exis_syst_equi}
	If $\phi_c(\cdot)$ and $\phi_g(\cdot)$ are strictly convex functions and $t_a(\cdot)$ is strictly monotone increasing, the system has a unique equilibrium $(\boldsymbol{c}^*,\boldsymbol{g}^*,\boldsymbol{\hat{x}}^*,\boldsymbol{\check{x}}^*, \boldsymbol{x}^*,  \boldsymbol{q}^*, \rho^{*})$ satisfying Definition \ref{defn:Syst_Equi}, where $(\boldsymbol{c}^*,\boldsymbol{g}^*,\boldsymbol{\hat{x}}^*,\boldsymbol{\check{x}}^*, \boldsymbol{x}^*,  \boldsymbol{q}^*)$ is the solution from the model~\eqref{eq:combined}, and $\rho_{\boldsymbol{\xi}}^k=\frac{\lambda_{\boldsymbol{\xi}}^k}{\pi_{\boldsymbol{\xi}}}$, for every $\boldsymbol{\xi}$ and $k$. 
\end{theorem} 

\state Proof. See \ref{app:pfs}. 

Convex reformulation \eqref{eq:combined} can be directly solved by commercial nonlinear solvers (e.g., IPOPT), which has provided an effective way of finding a system equilibrium as described in Definition~\ref{defn:Syst_Equi}. As the dimensions of uncertainties increase, the problem may become more challenging to solve. But since model \eqref{eq:combined} is convex, it can be solved by the classic scenario decomposition approach, such as the progressive-hedging (PH) algorithm.

\section{Numerical Examples}

\begin{figure}[htbp]
\begin{center}
    \includegraphics[width=1\textwidth]{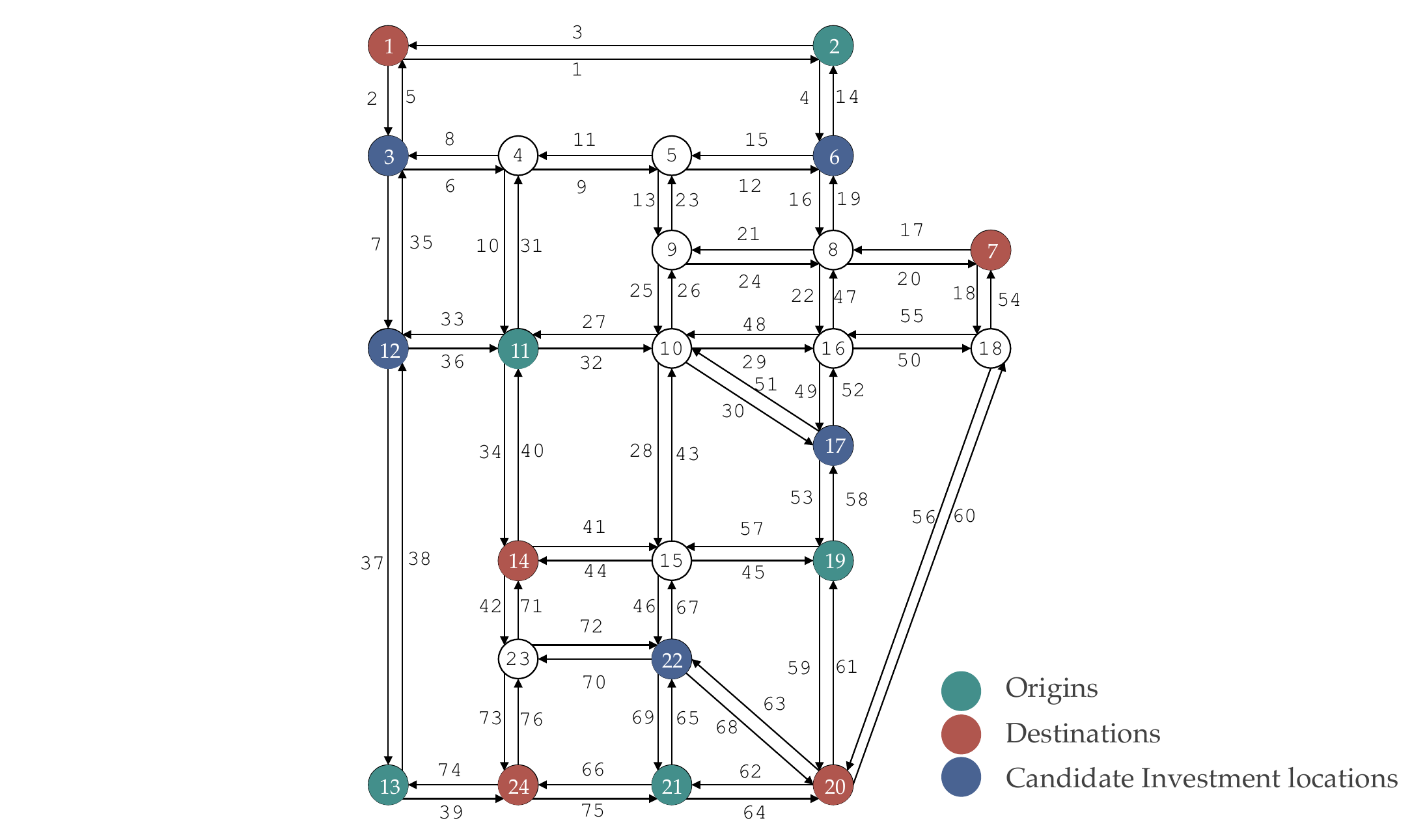}
\caption{Base Case Sioux Falls Test Network}
\label{fig:siou_fall}
\end{center}
\end{figure}

We use the Sioux Falls network, a widely used benchmark network, as shown in Figure \ref{fig:siou_fall}, to test the numerical performance of our solution methods and draw practical insights. The network consists of 24 nodes and 76 directed links. The number on each node/link in Figure \ref{fig:siou_fall} is the node/link index.

In Figure \ref{fig:siou_fall}, the green, red, and blue nodes (5 of each) represent the set of origins, destinations, and candidate facility locations, respectively. We consider 25 o-d pairs, each expecting 100 units of travel demand. For the link travel costs function, $t_a(v_a)$, we adopt a $4^{th}$-order Bureau of Public Roads (BPR) function: $t_a = t_a^0[1+0.15*(v_a/c_a)^4]$, where $t_a^0$ is the free flow travel time (FFT) and $c_a$ is the link capacity parameter\footnote{Note that $c_a$ is the ``capacity'' parameter used in BPR rather than the true link capacity}.  Values of $t_a^0$ and $c_a$ are documented in Table \ref{tab:capa_fft_base} in \ref{App:data_inputs}. For illustration purposes, the parameters in the travelers' utility function are assumed to be $\beta_0=0$, $\beta_1=1$, $\beta_2=0.06$, $e=1$. These utility parameter settings represent the case when users consider travel time and service prices when they choose facilities and routes and do not have a particular locational preference, and all users have a homogeneous service demand and value of time.  For facility providers, we select a quadratic form for the investment and operational cost function: $\phi_c(c) = 0.1c^2 + 170c$ and $\phi_g(g) = 0.1g^2 + 130g$. We refer to the above specifications as the base case, on which sensitivity analyses will be further conducted. Note that the magnitude of the parameters are arbitrary and just for illustration purposes. All the numerical experiments presented in this section were run on a 3.5 GHz Intel Core i5 processor with 8 GB of RAM under the Mac OS X operating system. 

\subsection{Deterministic Case}

\begin{figure*}[htbp]
	\centering
	\makebox[1.2\linewidth][c]{\hspace{-9em}
		\begin{subfigure}[t]{.46\linewidth}
			\centering
			\caption{Ignoring Congestion Case}
			\includegraphics[width=1\linewidth]{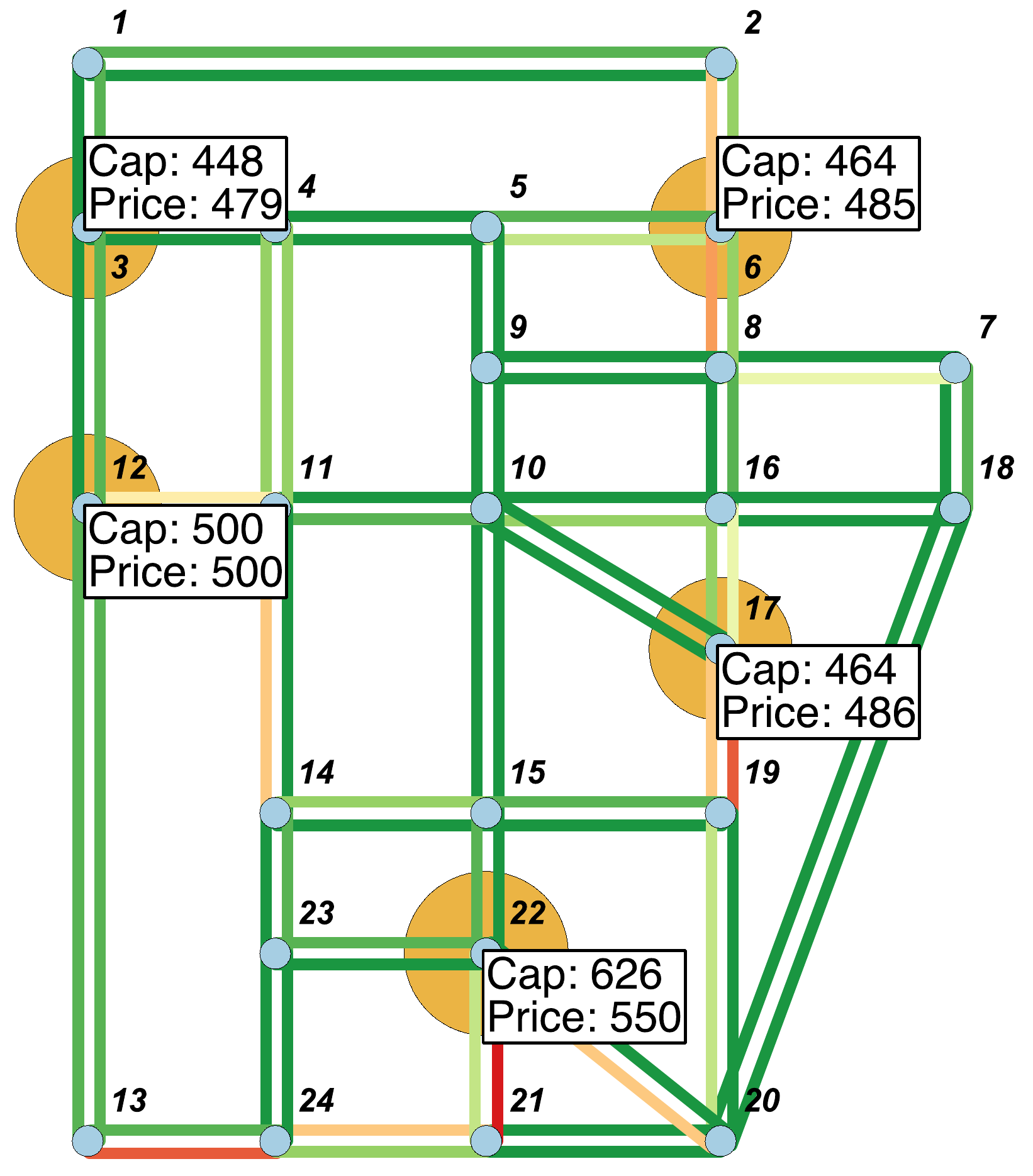} 
			\label{fig:no_cong_netw}
		\end{subfigure}
		\begin{subfigure}[t]{.28\linewidth}
			\centering
			\caption{System Evolution}
			\includegraphics[width=1\linewidth]{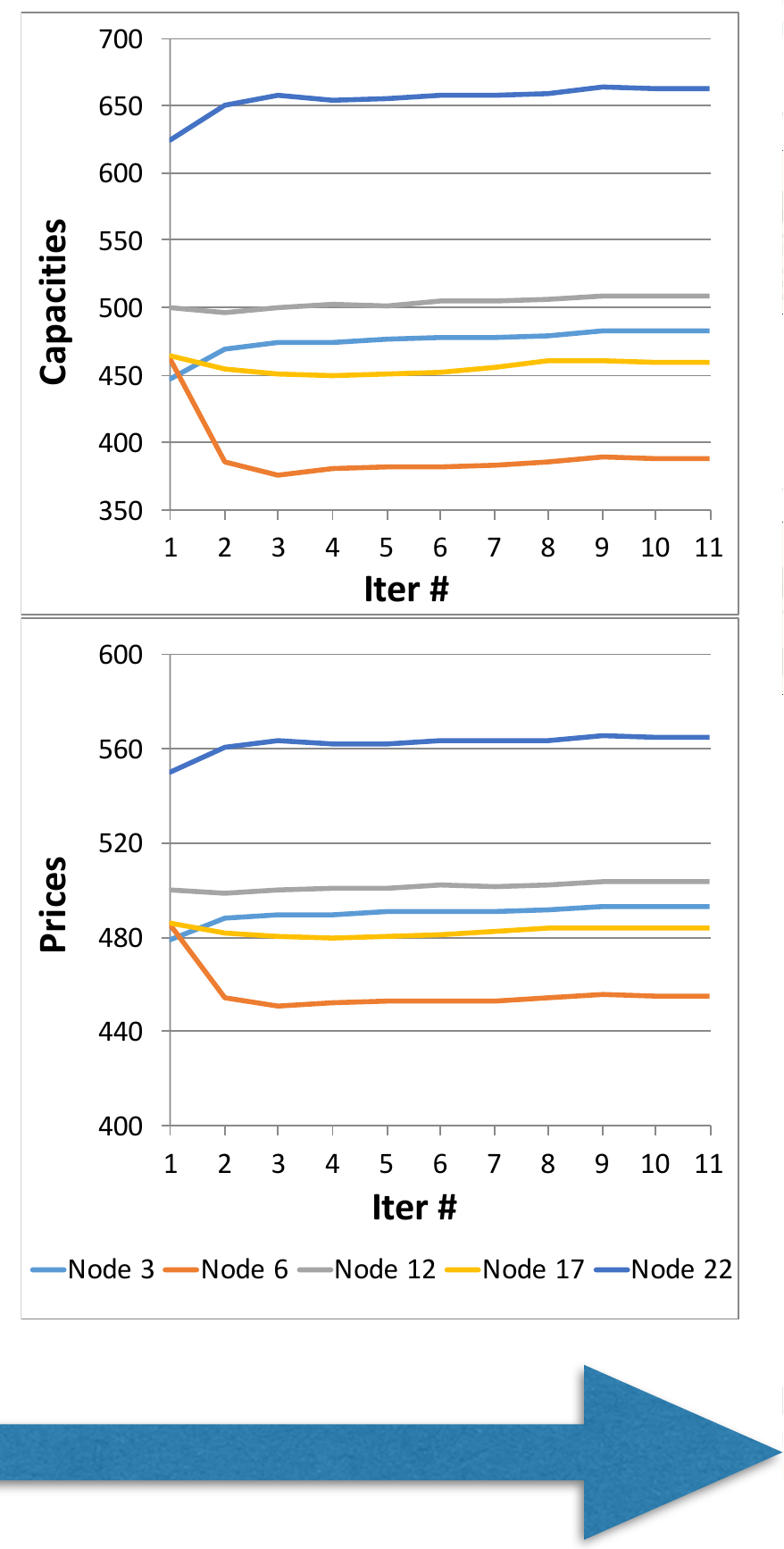} 
			\label{fig:cong_price_cap_conv}
		\end{subfigure}
		\begin{subfigure}[t]{.46\linewidth}
			\centering
			\caption{Base Case}
			\includegraphics[width=1\linewidth]{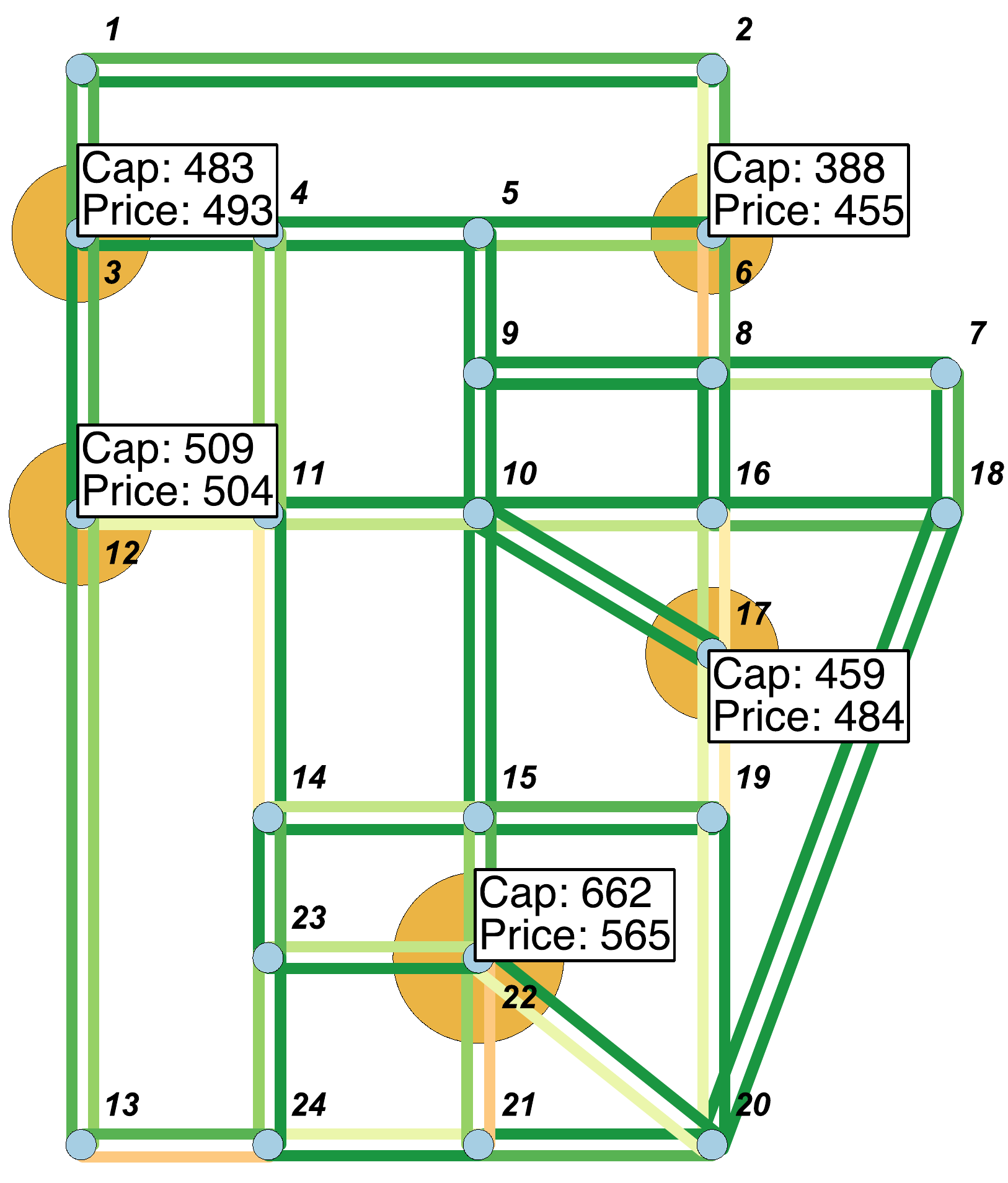} 
			\label{fig:base_case_netw}
	\end{subfigure}} \par\medskip
	\vspace{-2em}
	\makebox[1.0\linewidth][c]{\hspace{-0em}
		\begin{subfigure}{0.5\linewidth}
			\centering
			\includegraphics[width=0.7\linewidth]{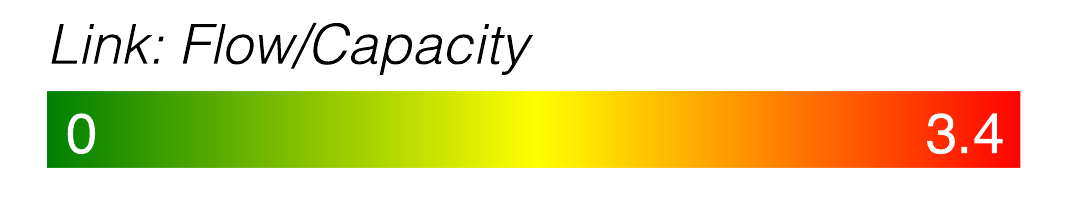} 
			\label{fig:netw_cong_lege}
	\end{subfigure}}
	\caption{Impacts of Modeling Network Congestion}
	\label{fig:impa_tran_cong}
\end{figure*}

First, we investigate the deterministic case where investors make an investment decision based on base-level future EV travel demand. \rev{Road congestion is an important factor in transportation system that can influence the decision making of investors as it will influence the facility selection and route choice of travelers \citep{duranton2011fundamental}. However, most of the existing studies do not consider traffic congestion in their planning of intermediate service facility (e.g., \citep{wang2019designing,xu2020mitigate}). To demonstrate quantitatively potential bias of investment resulted from lacking considerations of congestion, we compare numerical results between cases with and without congestion.} Figure \ref{fig:impa_tran_cong} illustrates the impacts of consideration of network congestion on model results, including the equilibrium capacity, price, and the traffic distribution of the transportation network.  Figure \ref{fig:no_cong_netw} corresponds to a case where transportation network congestion is not \rev{considered in the investors' decision making} (by setting the link capacity to be infinity), while Figure \ref{fig:base_case_netw} corresponds to the base case where network congestion is explicitly modeled. \rev{Note that, for the sake of fair comparison, the link capacity used for reporting the flow to capacity ratio in both figures are set to be equal to the actual capacity of the link (i.e., the links capacity in base case). The link flow in Figure \ref{fig:no_cong_netw} are hypothetical flow if congestion is not modeled in the intermediate facility location problem.}  We find that link congestion in Figure \ref{fig:base_case_netw} is not as significant as Figure \ref{fig:no_cong_netw} because when congestion is captured in the model, users will adjust their facility location and routing choices to avoid congestion. In addition, considering congestion costs for facility location choice will lead to a difference in the facility demand distribution, which leads to the difference in equilibrium service prices shown in Figure \ref{fig:impa_tran_cong}.  

Considering a hypothetical scenario, if investors planned the infrastructure without considering transportation network congestion (i.e., following the facility capacity results presented in Figure \ref{fig:no_cong_netw}), once the users experience significant congestion, their actual location and routing choices would be adapted, which means the computed equilibrium capacity and service prices would no longer be optimal for individual investors. Therefore, the investors would have an incentive to make changes accordingly.  If we let this system evolve, the system eventually converges to an equilibrium state identical to the base case results, as shown in Figure \ref{fig:cong_price_cap_conv}.  This experiment illustrates the importance of capturing realistic user-infrastructure interactions (e.g., user choices and transportation congestion) to find a stable system equilibrium state.

%

%
%

\begin{figure}[htbp]
\begin{center}
    \includegraphics[width=0.8\textwidth]{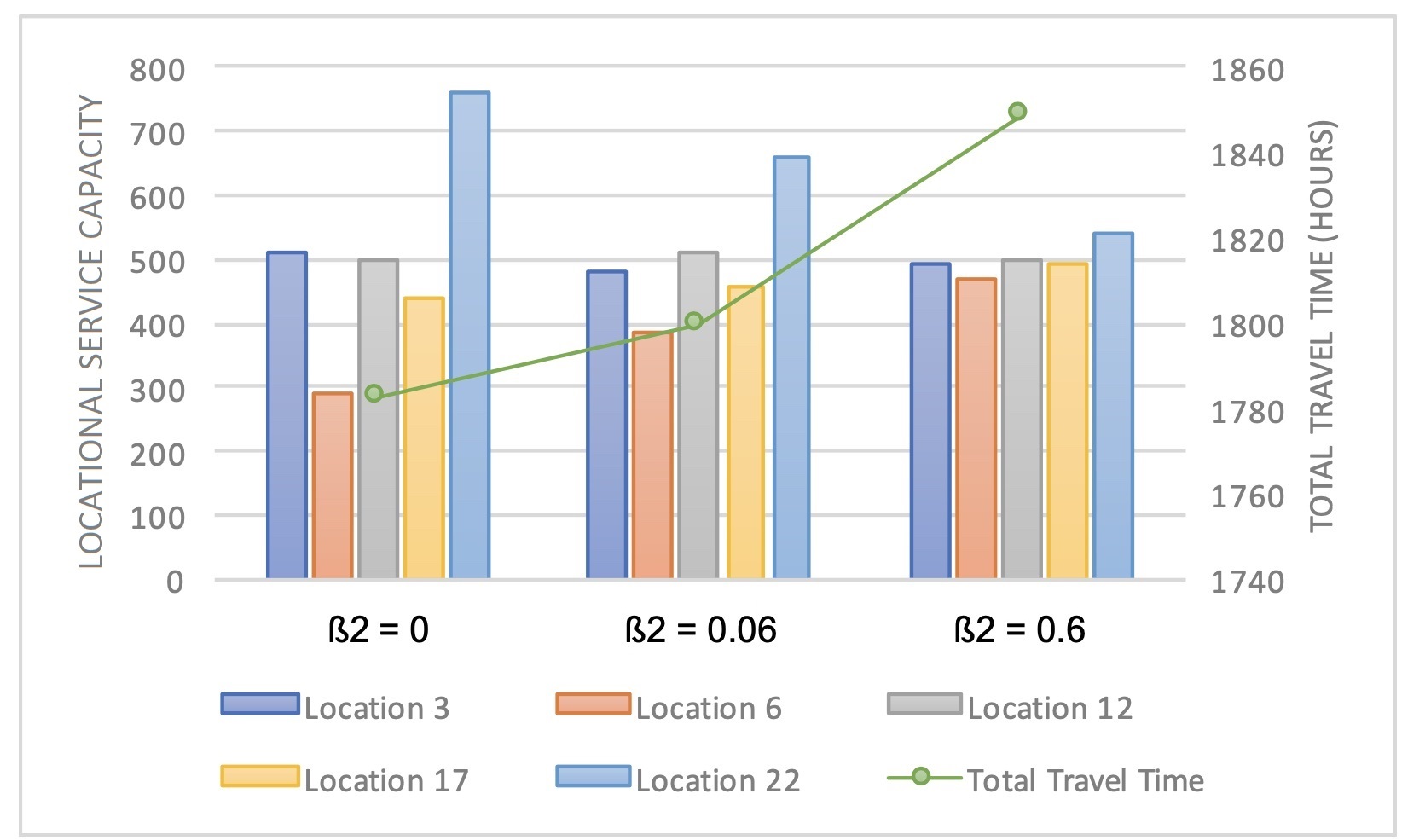}
\caption{\rev{Impact of User Preferences on Investment}}
\label{fig:pref_inve}
\end{center}
\end{figure}

Next, we study how different price sensitivities might affect the equilibrium investment layout and travel time by comparing cases with $\beta_2$ values at  0, 0.06, and 0.6. \rev{Figure \ref{fig:pref_inve} represents sensitivity analysis on $\beta_2$ and the specific magnitudes only aim to demonstrate the impact of different values of time on the equilibrium outcomes. More specifically, the interpretation of $\beta_1$ and $\beta_2$ are the dis-utility per unit of time and costs, respectively. Therefore, $\beta_1/\beta_2$ represents the value of time, i.e., monetary costs per unit of time.  Since $\beta_1 = 1$, $\beta_2 = 0$ represents the case when the value of time is infinitely, which means users choose the facility that takes the least detour and do not care about service prices. This case could represent facility providing life-critical service in emergency. When $\beta_2$ becomes larger, users put more weight on the service costs in addition to travel time when they choose facility. This case could represent daily facility service, such as EV charging facility.}  With increasing price sensitivity, the system has an increasing total travel time in equilibrium. The reason is that when travelers are more sensitive to price, they are more willing to choose a cheaper, albeit farther or more congested, service facility. In terms of equilibrium investment, higher price sensitivity leads to a more evenly distributed investment pattern because the preference of travelers for cheaper locations will naturally drive closer the equilibrium prices of different locations so that each location has similar attractiveness to investors.

\subsection{Stochastic Case}
\rev{Different sources of uncertainties can influence the choice of facility providers and users. The impact of uncertainty will be more prominent in the multi-agent framework \rev{because the response of each agent to the uncertainties will also influence the decisions of other agents.} \rev{In this section, we will focus on the uncertain total travel demand from each OD pair $d_{rs}$.} Other uncertainty sources can be similarly investigated.} 

The service demand is modeled with a random coefficient $\theta_\xi \in [\theta_\min,\theta_\max]$ for each scenario $\xi$ multiplied by the base case service demand. In the Sioux Falls test network, we considered $\theta_\min$ and $\theta_\max$ to be 1 and 1.2, respectively, and generated 20 service demand scenarios from a uniform distribution. We compared the results in three cases:
\begin{itemize}
\item{Case 1: Deterministic problem, where only the expected elastic demand is considered.}
\item{Case 2: Stochastic problem, where possible scenarios of demand and their associated probabilities are modeled.}
\item{Case 3: Wait-and-see problem, where all stakeholders have perfect forecast of the uncertainty parameters when they make investment decisions (i.e., relaxing the capacity variable to be scenario dependant without non-anticipitivity constraints).}  
\end{itemize}

These three cases allow us to investigate the decision makings of the stakeholders under different information availability scenarios (Section \ref{sec:decision}), as well as to quantify the stochastic programming metrics (e.g., the value of  stochastic solutions (VSS) and the expected value of perfect information (EVPI) \citep{birge2011introduction}) for each individual stakeholder (Section \ref{sec:metrics}).

The simulation time for the three cases were 0.153, 7.938, and 3.019 minutes respectively. The increased simulation time for cases 2 and 3, compared to case 1, indicates the additional computational burdens of stochastic modeling. Additionally, the computation time in case 3 is lower than case 2 because case 3 is basically running the deterministic problem (case 1) repeatedly for the number of scenarios (20 scenarios here).


\subsubsection{Results on the Decision Making of Stakeholders} \label{sec:decision}
\begin{figure*}
	\centering
		\begin{subfigure}[t]{.6\linewidth}
			\centering
			\caption{}
			\includegraphics[width=1\linewidth]{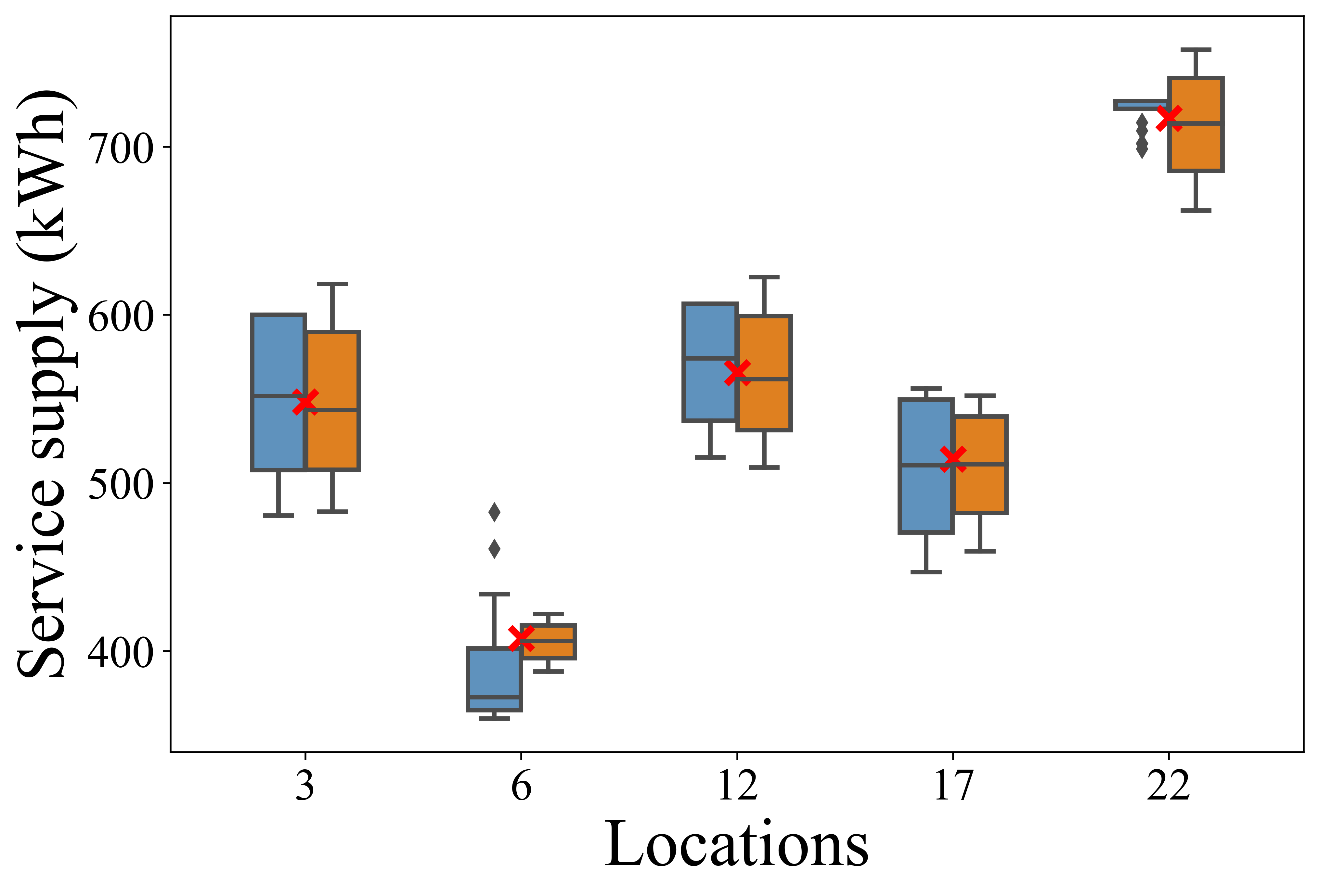} 
			\label{fig:result_services}
		\end{subfigure}\\
		
		\begin{subfigure}[t]{.6\linewidth}
			\centering
			\caption{}
			\includegraphics[width=1\linewidth]{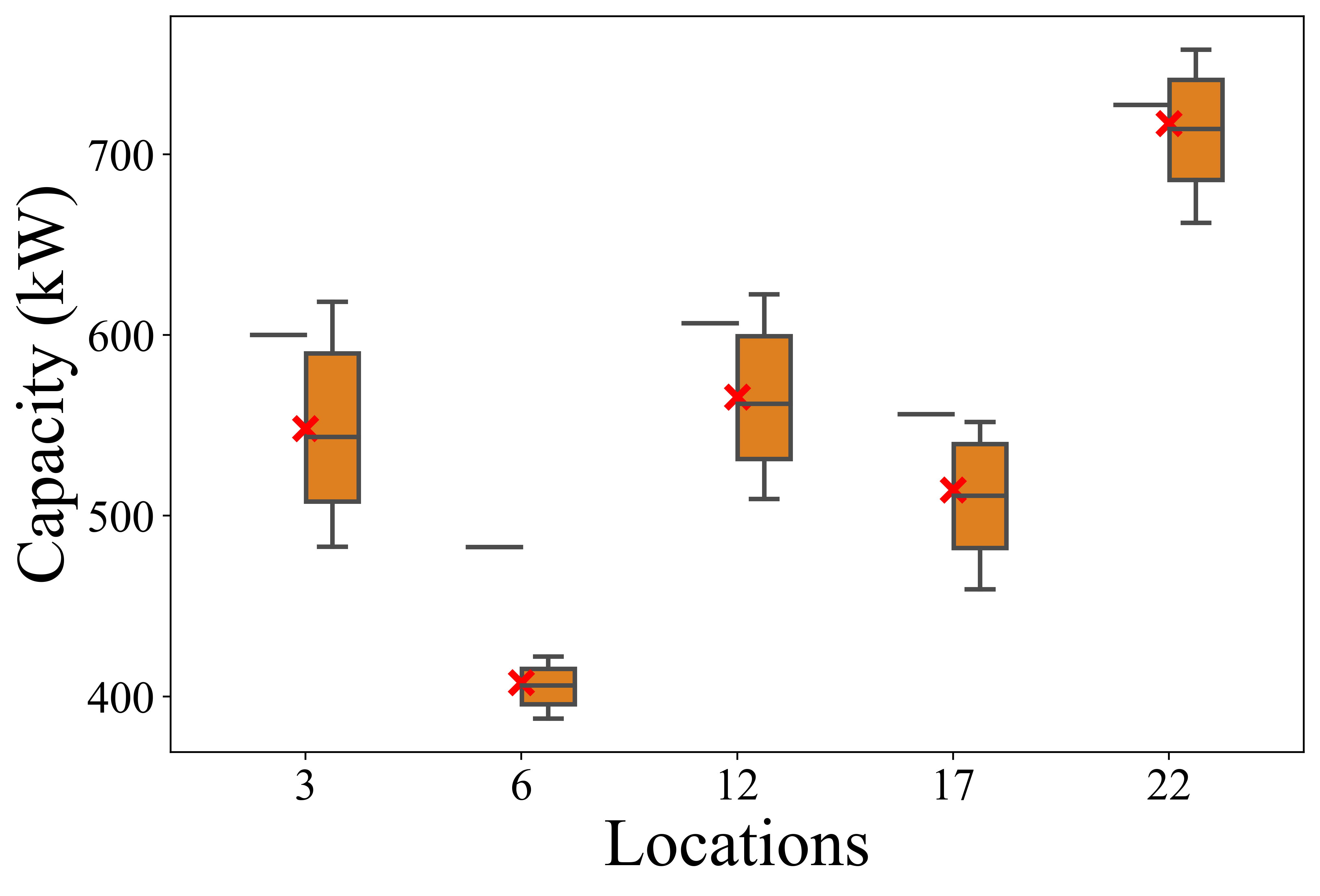} 
			\label{fig:results_capacity}
		\end{subfigure}
		\\
		
		\begin{subfigure}[t]{.6\linewidth}
			\centering
			\caption{}
			\includegraphics[width=1\linewidth]{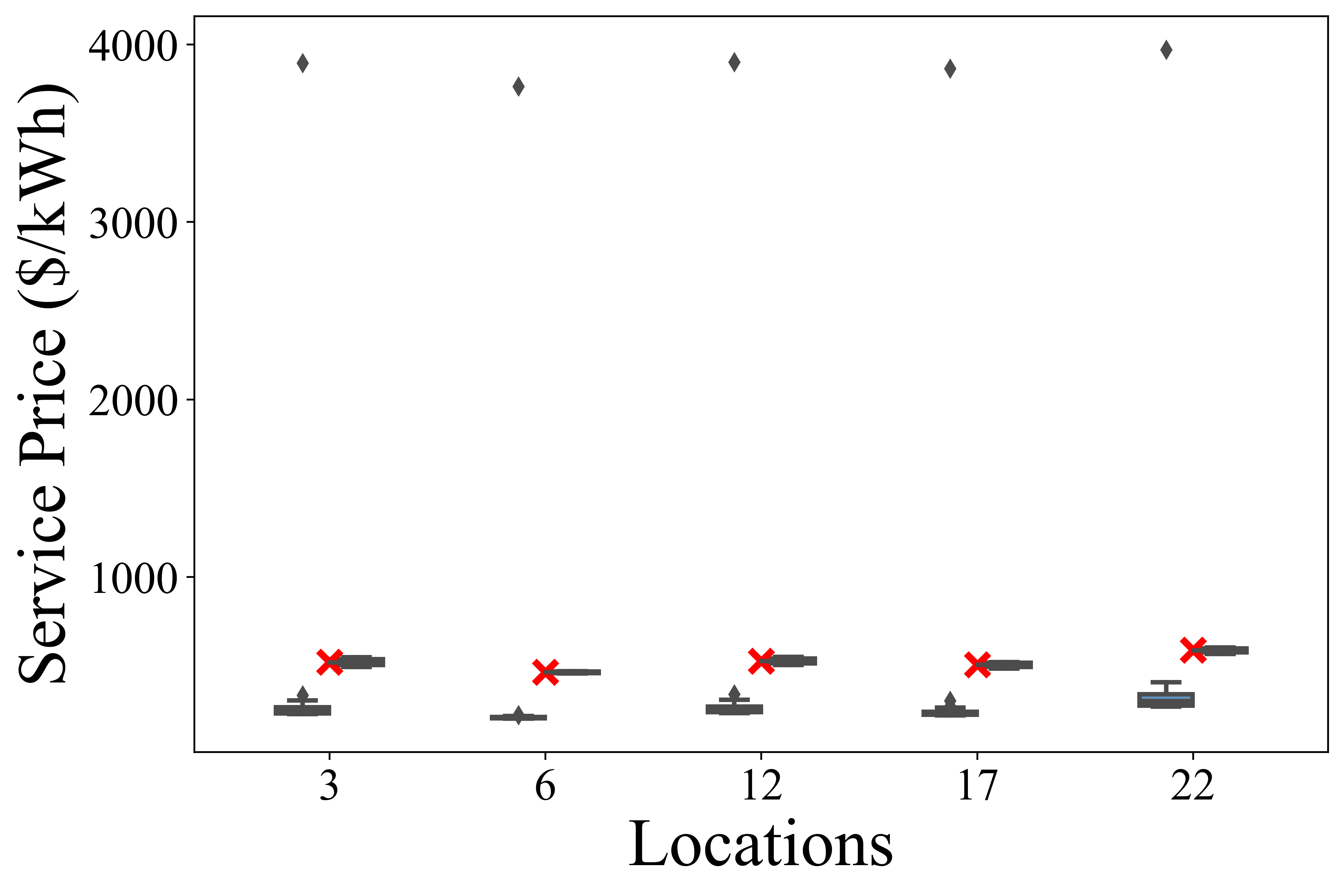} 
			\label{fig:results_prices}
	\end{subfigure}
	\vspace{-2em}
	\caption{Resulting decision variables (a) services provided by the facility locations, (b) determined capacity of facility locations, and (c) prices of the provided services;\textcolor{red}{{$\boldsymbol{\times}$}} case 1,  \Rectangle[Abi] case2, and \Rectangle[orange] case 3. }
	\label{fig:stochastic_results}
\end{figure*}

\rev{Supplied service quantities, service prices, and total capacity at facility locations for each of the three cases are presented in Figure \ref{fig:stochastic_results}.} From Figure \ref{fig:result_services}, we can see that the service supply quantities at all facility locations are similar between case 1 and the mean value in case 3 \rev{(i.e., ``the supply for mean demand is similar to the mean supply for each demand sceanrio'').}  However, the ranges of service supply in cases 2 and 3 vary across different locations. \rev{The main reason stems from agents being able to make different investment decisions for each scenario in case 3 and have more flexibility to optimize their supply quantity.} For example, services provided at location 6 have a higher supply on average in case 3 than in case 2, meaning that perfect information on OD demand can encourage facility provision in location 6. In terms of the variance of service supply, perfect information (case 3) does not have universal impacts for different locations. For example, the variance of service supply in case 3 is lower at location 6, while it is higher at location 22, compared with case 2. Smaller variance indicates that service supply is less sensitive to uncertain parameters. 



Similar to the supplied services, the locational capacities in case 1 are close to the mean capacity of the stochastic problem in case 3 (see Figure \ref{fig:results_capacity}), which is as expected since investors will invest the same amount of capacity as the supply quantity for each scenario in both case 1 and case 3 and there is no unused capacity. Enforcing the non-anticipativity constraints in case 2 results in a constant locational capacity for all the scenarios (see Figure \ref{fig:result_services}). \rev{In general, the capacity at each location in case 2 is higher than the mean of the capacity investment in case 3, because investors in case 2 consider all scenarios and may need to over-invest to consider the most profitable demand scenario, which has the highest service prices, to maximize their expected profits. In case 3, however, the investors determine capacities for each scenario independently and can invest less for scenarios where lower capacity is needed. Additionally, with perfect information, investors have the flexibility to distribute their capacity among locations for each scenario. For example, facilities at 3, 12, and 22 may have higher capacity investment in some extreme scenarios compared to case 2.} 

Figure \ref{fig:results_prices} shows the equilibrium service prices for each case. In general, the service prices in cases 1 and 3 are higher than the prices in case 2 (see Figure \ref{fig:results_prices}). However, in scenarios where the facility capacity is binding ($g^k_\xi=\bar{c}^k_\bullet$), we observe drastically higher service prices in case 2. These spikes in service prices stem from the market-based modeling framework, where the service prices are based on the marginal costs of additional service quantity. In case 2, the costs of the additional service quantity in the capacity-binding scenario will need to account for the costs of unused capacity in the other scenarios. The pricing mechanism based on marginal costs also explains why facility investors tend to invest more facing future uncertainties since they can make more profits during the supply shortage.

{In summary, the key observations of the case analyses are as follows:
\begin{itemize}
    \item \rev{With a perfect forecast on the uncertain parameters (case 3), the decision-makers act similarly on average in terms of supply capacity and quantity compared with the deterministic case  (case 1), as presented in Figures \ref{fig:result_services} and \ref{fig:results_capacity}. In other words,  the mapping from scenario to facility supply has the property that the supply for mean demand is similar to the mean supply for each demand scenario}.
    
    \item \rev{Stochastic decision-making (case 2) leads to more investment of facility capacity as a response to the uncertain future demand, which will lead to different locational service supplies and demand quantities compared with the case with perfect information on the uncertainties (case 3).} For example, Figure \ref{fig:results_capacity} shows significantly higher invested capacity (especially at location 6) in case 2 compared to case 3, which has resulted in different service supply quantities, as presented in Figure \ref{fig:result_services}.
    
    \item The market-based mechanism of service pricing could result in high prices for scenarios where the capacity constraints are binding in the stochastic problem (case 2), which leads to higher investment in case 2 compared with case 3. This is evident from the marginal locational service prices presented in Figure \ref{fig:results_prices}.
\end{itemize}
}

\subsubsection{Stochastic Programming Metrics} \label{sec:metrics}

In order to quantify the impacts of stochastic programming and uncertainty information on the benefits of individual stakeholders, we investigate two classic stochastic programming metrics in the context of network equilibrium: (1) value of stochastic solutions (VSS) and (2) expected value of perfect information (EVPI). 

VSS evaluates the potential benefit of implementing stochastic programming solutions considering system uncertainties compared with the deterministic solutions. EVPI is another stochastic metric that quantifies the value of the perfect forecast of the uncertain parameters on the decision-making of stakeholders. In order to calculate these metrics, we first calculate the objective value of each stakeholder with the results found from the previously defined three cases. Then, the VSS would be the objective value difference of each stakeholder between case 2 (stochastic problem) and case 1 (deterministic problem). Since case 3 models the condition where all stakeholders have access to the information of uncertain scenarios, EVPI would be the objective value difference between case 3 (the wait and see problem) and case 2. Note that in contrast to single-agent stochastic programming, VSS and EVPI may be negative in a multi-agent setting due to the complex interactions.

The objective value of facility providers can be determined using model (\ref{eq:investor}), and the objective of the service users will be calculated as the total expected utility based on equation (\ref{eq:total_utility}), where $U^{\mathrm{total}}$ is the total utility of the users, $\boldsymbol{U}$ is the drivers' utility defined in equation (\ref{eq:utility}), and $\boldsymbol{q}$ are the resulting traffic flow after solving the problem. 
\begin{align}
  U^{\mathrm{total}} =  \mathbb{E}_{\boldsymbol{\xi}}\sum_{k \in K} \frac{1}{\beta_2} ( q^{rsk}_\xi \times U^{rsk}_\xi)  \label{eq:total_utility}
\end{align}
Notice that in equation \eqref{eq:total_utility}, we have divided the expected utility by $\beta_2$ to normalize the utility in the unit of \$. Therefore, we can also analyze the system welfare or surplus as the summation of the service providers' objective and users' total utility.

\begin{figure*}
	\centering
			\includegraphics[width=0.5\linewidth]{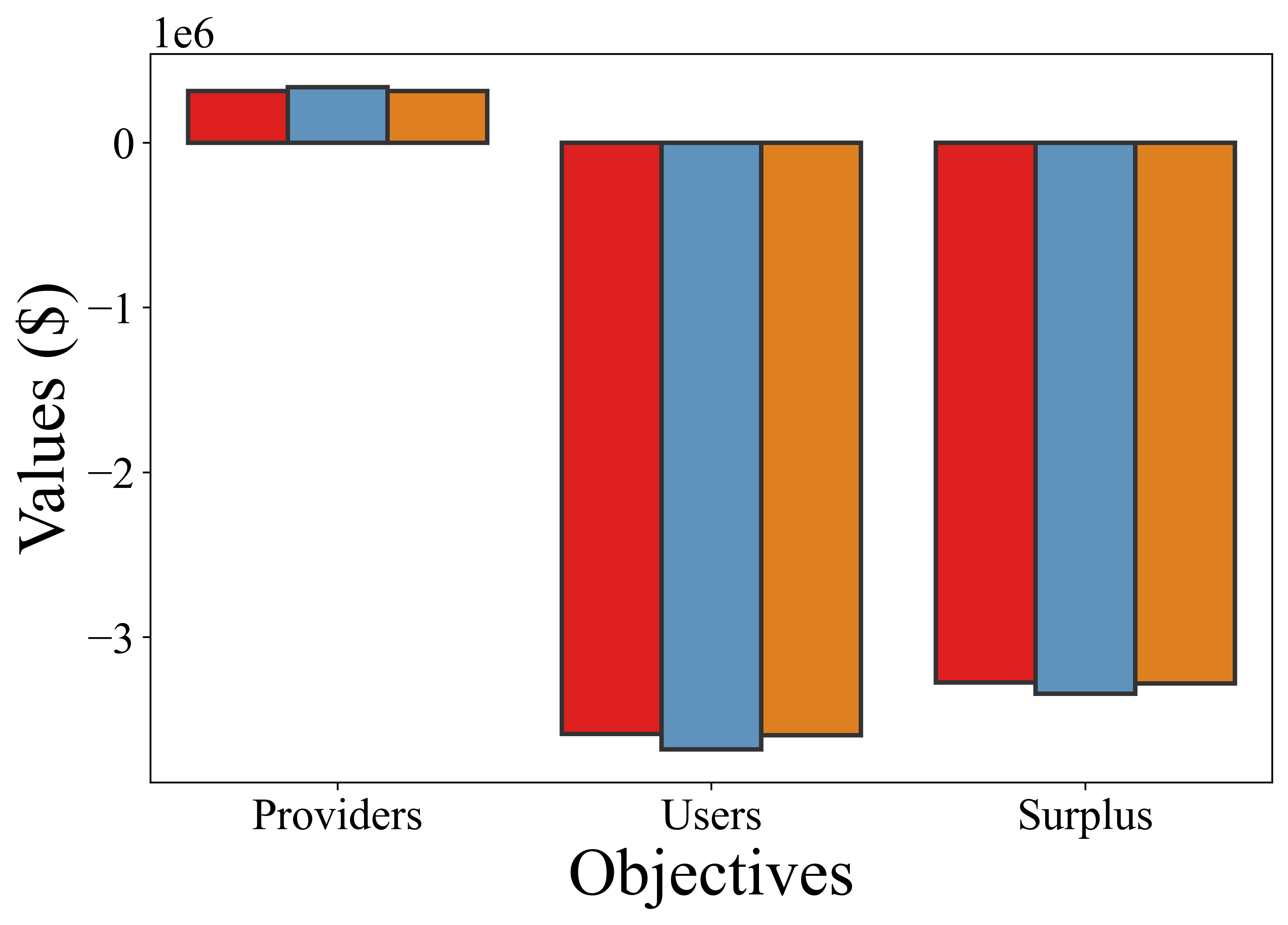} 
	\caption{Stakeholders' objectives and system surplus; \Rectangle[Red] case 1,  \Rectangle[Abi] case2, and \Rectangle[orange] case 3. }
	\label{fig:VPI}
\end{figure*}
The calculated objectives and system surplus are shown in Figure \ref{fig:VPI} for the three cases. \rev{We can see that facility providers have the highest objective value in case 2 while facility users have the lowest utility, because in case 2, the market will yield a much higher  service price  on facility users in capacity binding scenarios (Figure \ref{fig:results_prices}), resulting in higher benefits for investors and lower utility for the users.} The objective value of providers has decreased with perfect information (case 3) compared with case 2 mainly because of the high spikes of service prices in case 2 when the supply quantities are constrained by the capacity of the facilities. The objective value of users has become less negative, which means that users have benefited from perfect information. The improvement of objective value for users was more prominent than the providers resulting in improved system surplus in case 3 compared to case 2.

The mentioned stochastic metrics will help us solidify these comparisons. Based on the definition provided for VSS, the VSS for the providers would be 24652.7 units, which represents the profit of investors by implementing stochastic solutions compared to deterministic modeling. Comparing the objectives in cases 2 and 3, the objective value of providers has decreased by 23394.0 units, and the objective value of users has improved by 84225.6 units. These changes can be interpreted as EVPI. The metrics also show that providers are better off without having perfect information over the realization of uncertain parameters, whereas users have benefited from the perfect information. From a system perspective, perfect information (case 3) and deterministic solutions (case 1) achieve higher system surpluses compared with the stochastic solution (case 2).

\subsection{ Results for disaggregated investors}
    
\rev{The results presented in the previous section for investment capacity is based on the aggregated cost function of investors at each facility location. Here, we will explicitly model multiple investors with heterogeneous cost functions to investigate the equilibrium capacity share between investors. The objective function of each investor will take the same quadratic forms}
\begin{align}
    \Phi(c_{i}^k) = a_i {c_{i}^k}^2 + b_i {c_{i}^k}, 
\end{align}
\rev{with $a_i$ and $b_i$ being the cost coefficients for each investor $i$. 

As an example, we will consider having two disaggregated investors $i\in \{1,2\}$. The first investor with cost coefficients $a_1$ = 0.1 and $b_1$ = 170, and the second investor with $a_2$ = 1 and $b_2$ = 17. The first investor represents an investor who has higher marginal costs for installation but the marginal costs increase slower with capacity, and the second investor represents an investor with the opposite cost structure. An example for the first investor would be a firm that owns large area of land with no existing infrastructure and an example for the second investor would be a firm that owns limited land with underlying infrastructure already installed (e.g., gas station owners). 

Figure \ref{fig:Inv_share} shows the capacity share between these two investors for the stochastic case (case 2). Since, the second investor has higher investment cost for high capacity, it has lower share of investment capacity at each location. The differences between total capacity installed among different locations depend on the route and path selection of drivers as we have discussed in the stochastic and deterministic analysis. Note that we can use an aggregated cost function to model the total capacity at each node. This is evident in our example that the total capacity at each location presented in Figure \ref{fig:Inv_share} is equal to the total capacity based on the aggregated cost function.}

\setcounter{figure}{7}    
\begin{figure}[H]
        \begin{center}
        \includegraphics[width=0.6\linewidth]{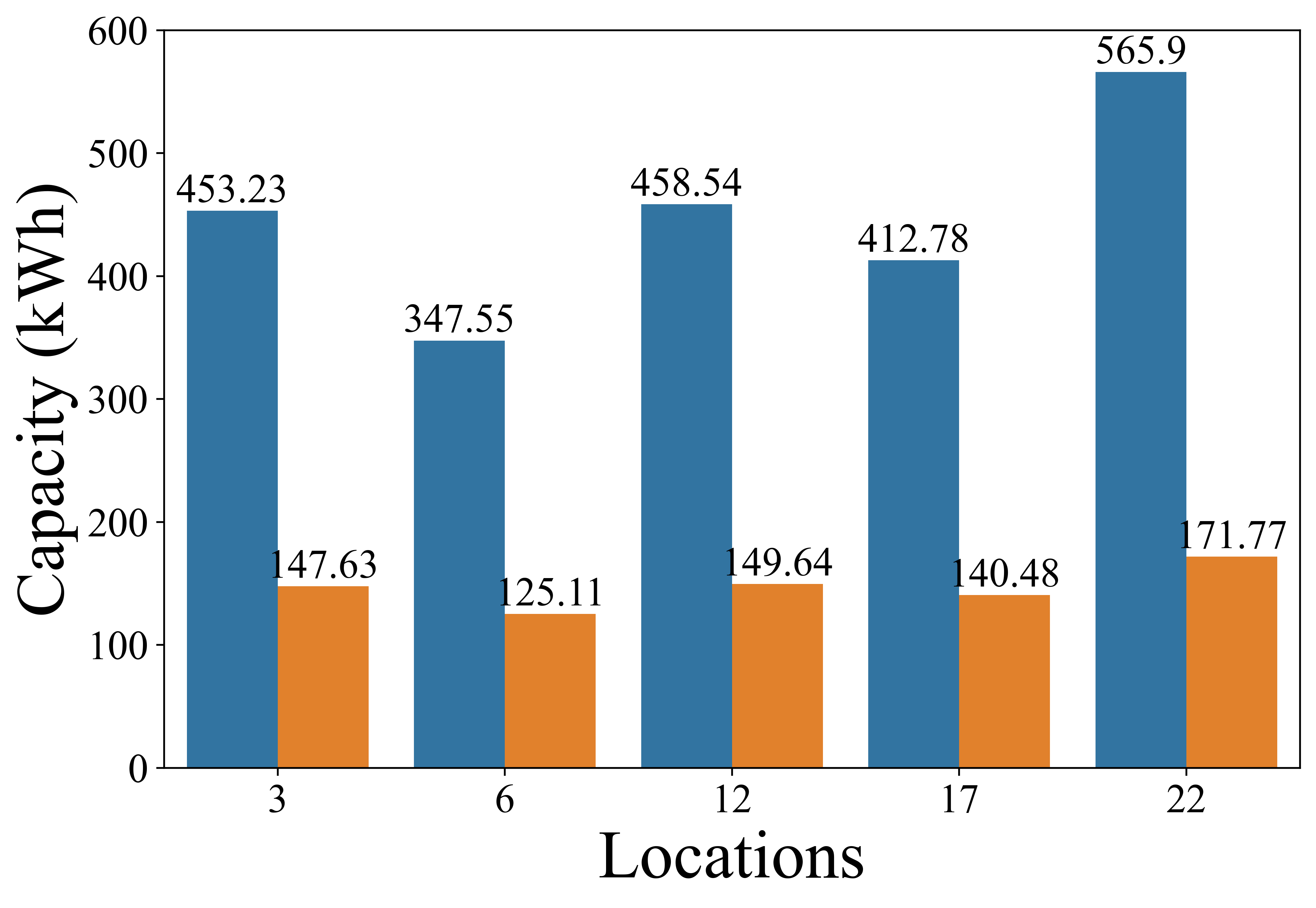}
        \caption{\rev{Investors' capacity shares in facility locations.}}
        \label{fig:Inv_share}
        \end{center}
\end{figure}

\section{Discussion}

In this paper, we have presented a new modeling framework along with an efficient solution method for long-term infrastructure system planning problems with challenges brought by intermediate facility provision, non-cooperative stakeholders, and complex agent-infrastructure interactions. The existence and uniqueness of equilibrium are proved. Through numerical examples, we demonstrated that (1) prices, investment, and profits might differ significantly across locations due to agent-infrastructure interactions; (2) ignoring transportation congestion may lead to system assessment bias; (3) the equilibrium investment patterns may be sensitive to user preferences; (4) while the stochastic decision making facing uncertainties may lead to higher investment compared with the average investment with perfect information, it will lead to significantly higher service prices in the scenarios when the capacity constraint is binding; and (5) information on uncertain parameters may not benefit the facility providers, but will benefit the users and increase the system surplus.

This research can be extended in several directions.  From a methodological viewpoint, the assumption of perfect competition may not fit all the applications when there is noticeable market power from the supply or demand side. Therefore, different market structures could be investigated.  For example, if suppliers have the strong market power to anticipate and influence the responses of users, a bilevel leader-and-follower model would be more appropriate. In addition, we can leverage the existing modeling framework to investigate the control strategies to influence market interactions so that a more resilient facility network can be achieved. Also, we observed asymmetry of value of information for the service providers versus service users, which is unique in a decentralized decision environment.  How to design an effective information acquisition and sharing mechanism in this context would be an interesting topic for future investigation. 

\section*{References}
\bibliography{all.bib}
\pagebreak 
\appendix
\section{Proofs.}\label{app:pfs}

\state Proof (Lemma \ref{lem:GCDA}).
Firstly, the objective function (\ref{obj:CDA}) is convex becasue it is a linear combination of three basic convex functions: (1) $f_1(x) =  \int_{0}^{x} g(u) \mathrm{d}u$, with $g(u)$ being a positive and nondecreasing function, (2) $f_2(x) = x\ln x$ and (3) $f_3(x) = cx$. In addition, the constraints for problem (\ref{eq:cda}) are all linear. Therefore, the optimization problem (\ref{eq:cda}) is convex. Because of the differentiability of function (\ref{obj:CDA}), the optimality conditions of problem (\ref{eq:cda}) is equivalent to the following complementarity conditions in additions to constraints (\ref{cons:x_p} $\sim$ \ref{cons:q_d}): 
$\forall a\in \mathcal{A}, r \in R, s \in S, k \in K^{rs}, p \in P^{rsk}$
\begin{subequations}
\begin{align}
 0\leq x_p& && \perp && \sum_{a \in \mathcal{A}_p} t_a(\cdot) - \boldsymbol{\gamma}^T(B_{\hat{p}}+B_{\check{p}}) \geq 0 \label{kkt_x_p}\\
 0\leq \hat{x}_a^{rsk}& &&\perp && \gamma_a - A_a^T\boldsymbol{\hat{\lambda}} \geq 0 \label{kkt_x1_a}\\
 0\leq \check{x}_a^{rsk}& &&\perp && \gamma_a - A_a^T\boldsymbol{\check{\lambda}} \geq 0 \label{kkt_x2_a}\\
 0\leq q^{rsk}& &&\perp &&\frac{1}{\beta_1}(ln(q^{rsk}) + \beta_3\frac{\rho^ke^{rs}}{inc^{rs}} - \beta_2\sum_{i\in I_k}c_i^s - \beta_0^k) + E^{rk T}\boldsymbol{\hat{\lambda}} + E^{ks T}\boldsymbol{\check{\lambda}} + \mu^{rs} \geq 0\label{kkt_q}
\end{align}
\label{kkt}
\end{subequations}

We first show that the traffic flow solutions is Wardrop user equilibrium by proving the following two conditions. 

\begin{description}
\item[1. All the used paths connecting $r, s, k$ have the same travel time. ] $\forall r \in R, s \in S, k \in K^{rs}$, for those $\tilde{p} \in P^{rsk}$ with $x_{\tilde{p}} > 0$,  $\sum_{a \in \mathcal{A}_{\tilde{p}}} t_a(\cdot) = \boldsymbol{\gamma}^T(B_{\hat{{\tilde{p}}}}+B_{\check{{\tilde{p}}}})$ (because of  (\ref{kkt_x_p})). Due to the following two conditions: 
\begin{itemize}
\item  for $\tilde{a} \in \hat{\tilde{p}}$, i.e. $B_{\hat{{\tilde{p}}},\tilde{a}} = 1$, $\hat{x}_{\tilde{a}}^{rsk} > 0$ and therefore $\gamma_{\tilde{a}} = A_{\tilde{a}}^T\boldsymbol{\hat{\lambda}}$ (because of (\ref{kkt_x1_a})). So $\boldsymbol{\gamma}_{\tilde{a}}^TB_{\hat{{\tilde{p}}},{\tilde{a}}} = \boldsymbol{\hat{\lambda}}^TA_{{\tilde{a}}}B_{\hat{{\tilde{p}}},{\tilde{a}}}$
\item for $\tilde{a} \notin \hat{\tilde{p}}$, i.e. $B_{\hat{{\tilde{p}}},\tilde{a}} = 0$, $\boldsymbol{\gamma}_{\tilde{a}}^TB_{\hat{{\tilde{p}}},{\tilde{a}}} = \boldsymbol{\hat{\lambda}}^TA_{{\tilde{a}}}B_{\hat{{\tilde{p}}},{\tilde{a}}} = 0$
\end{itemize}
we have $\boldsymbol{\gamma}^TB_{\hat{\tilde{p}}} = \boldsymbol{\hat{\lambda}}^TAB_{\hat{\tilde{p}}} $. Notice that $AB_{\hat{p}} = E^{rk} $, so $\boldsymbol{\gamma}^TB_{\hat{\tilde{p}}} = \boldsymbol{\hat{\lambda}}^TE^{rk} $. 

Same procedure, we have $\boldsymbol{\gamma}^TB_{\check{\tilde{p}}} = \boldsymbol{\check{\lambda}}^TE^{ks} $. 

So $\sum_{a \in \mathcal{A}_{\tilde{p}}} t_a(\cdot) = \boldsymbol{\gamma}^T(B_{\hat{{\tilde{p}}}}+B_{\check{{\tilde{p}}}}) =\boldsymbol{\hat{\lambda}}^TE^{rk} + \boldsymbol{\check{\lambda}}^TE^{ks} \doteq \tau^{rsk} $, which only dependents on $r,s,k$.
\item[2. All the unused paths connecting $r, s, k$ have no smaller travel time]
{\bf than that of the used paths.} $\forall r \in R, s \in S, k \in K^{rs}$, for those $\tilde{p} \in P^{rsk}$ with $x_{\tilde{p}} = 0$,  $\sum_{a \in \mathcal{A}_{\tilde{p}}} t_a(\cdot) \geq \boldsymbol{\gamma}^T(B_{\hat{{\tilde{p}}}}+B_{\check{{\tilde{p}}}})$ (because of  (\ref{kkt_x_p})). From (\ref{kkt_x1_a}, \ref{kkt_x2_a}), $\gamma_{\tilde{a}} \geq A_{\tilde{a}}^T\boldsymbol{\hat{\lambda}}$ and  $\gamma_{a} \geq A_{a}^T\boldsymbol{\check{\lambda}}, \forall a$. So $\sum_{a \in \mathcal{A}_{\tilde{p}}} t_a(\cdot) \geq \boldsymbol{\gamma}^T(B_{\hat{{\tilde{p}}}}+B_{\check{{\tilde{p}}}}) \geq \boldsymbol{\hat{\lambda}}^TAB_{\hat{\tilde{p}}} + \boldsymbol{\check{\lambda}}^TAB_{\check{\tilde{p}}}=\boldsymbol{\hat{\lambda}}^TE^{rk} + \boldsymbol{\check{\lambda}}^TE^{ks} = \tau^{rsk}$.
\end{description}

Next, we show the OD-demand solutions are the service location choice with logit facility demand functions. This can be easily seen from (\ref{kkt_q}): for any $k$ with $q^{rsk} > 0$, $\frac{1}{\beta_1}(ln(q^{rsk}) + \beta_3\frac{\rho^ke^{rs}}{inc^{rs}} - \beta_2\sum_{i\in I_k}c_i^s - \beta_0^k) + E^{rk T}\boldsymbol{\hat{\lambda}} + E^{ks T}\boldsymbol{\check{\lambda}} + \mu^{rs} =0$. After reorganization,
\begin{align*}
q^{rsk} &=e^{\beta_0^k - \beta_1(E^{rk T}\boldsymbol{\hat{\lambda}} + E^{ks T}\boldsymbol{\check{\lambda}}) +  \beta_2\sum_{i\in I_k}c_i^s - \beta_3\frac{\rho^ke^{rs}}{inc^{rs}} + \beta_1\mu^{rs}}\\
&=e^{\beta_0^k - \beta_1\tau^{rsk} +  \beta_2\sum_{i\in I_k}c_i^s - \beta_3\frac{\rho^ke^{rs}}{inc^{rs}} + \beta_1\mu^{rs}}\\
&=e^{U^{rsk} + \beta_1\mu^{rs}}
\end{align*}
\eop

\state Proof of Lemma~\ref{lem:convex_existence}. 

Objective function \eqref{obj:centralized} is a linear combination of five types of function: $\phi_c(\cdot)$, $\phi_g(\cdot)$, $\int_{0}^{v_{a,\boldsymbol{\xi}}}t_{a}u \textrm{d}u$, $q\textrm{ln} q$, and $q$. First, the representative investor costs functions, $\phi_c(\cdot)$ and $\phi_g(\cdot)$ are convex by assumption. Second, since the link performance function $t_a(\cdot)$ is monotone increasing, thus $\int_{0}^{v_{a,\boldsymbol{\xi}}t_{a}u du}$ is a convex function because its second order derivative $t_a^{'}(\cdot) \geq 0$. Third, the rest of the functions $q\ln q$ and $q$ can be easily shown to be convex by taking second order derivative.  Therefore, the reformulation problem (model \eqref{eq:combined}) corresponds to a minimization of a convex function under linear constraints. Therefore, model \eqref{eq:combined} is a convex problem \citep{VaAn}. Under constraint qualifications, the convex reformulated problem (i.e., model \eqref{eq:combined}) has at least one solution. 

Furthermore, if $\phi_c(\cdot)$, $\phi_g(\cdot)$ are strictly convex functions and $t_a(\cdot)$ is strictly monotone increasing, following the same logic as above, model \eqref{eq:combined}) corresponds to a minimization of a strictly convex function under linear constraints. Therefore, the solution of  model \eqref{eq:combined}), if exists, is unique.

\eop

\state Proof of Lemma~\ref{lem:equivalency}. 

  Since $\phi_c(\cdot)$ and $\phi_g(\cdot)$ are convex functions and $t_a(\cdot)$ is monotone increasing, following Lemma~\ref{lem:convex_existence}, model \eqref{eq:combined} is a convex optimization problem. It is not difficult to see that the first order conditions associated with the convex optimization problem~\eqref{eq:combined} with dual multipliers $\{\lambda_\xi^k\}$ are separable by agents, and each corresponds to the first order conditions of the representative investor problem~\eqref{obj:inve}-\eqref{cons:capa} and the GCDA problem~\eqref{obj:CDA}-\eqref{cons:q_d} respectively, with equilibrium price vector $\rho_{\boldsymbol{\xi}}^k=\frac{\lambda_{\boldsymbol{\xi}}^k}{\pi_{\boldsymbol{\xi}}}$, for every $\boldsymbol{\xi}$ and $k$, where $\{\pi_{\boldsymbol{\xi}}:\boldsymbol{\xi} \in \boldsymbol{\Xi}\}$ is the probability distribution of $\boldsymbol{\xi}$.
\eop

\state Proof of Theorem~\ref{thm:exis_syst_equi}. 

  Theorem~\ref{thm:exis_syst_equi} directly follows from Lemma~\ref{lem:convex_existence} and Lemma~\ref{lem:equivalency}.
\eop

\pagebreak 
\setcounter{table}{0}
\section{Data inputs}\label{App:data_inputs}

\begin{table}[htbp]
  \centering
  \caption{Base Case Link Capacity $c_a$ (veh/h) and Free-flow Travel Time, FTT $t_a^0$ (min)}
  \resizebox{\columnwidth}{!}{
    \begin{tabular}{ccccccccccccccc}
    \toprule
    Link  & FFT   & Capacity  & \multicolumn{1}{c}{Link} & \multicolumn{1}{c}{FFT} & \multicolumn{1}{c}{Capacity } & \multicolumn{1}{c}{Link} & \multicolumn{1}{c}{FFT} & \multicolumn{1}{c}{Capacity } & \multicolumn{1}{c}{Link} & \multicolumn{1}{c}{FFT} & \multicolumn{1}{c}{Capacity } & \multicolumn{1}{c}{Link} & \multicolumn{1}{c}{FFT} & \multicolumn{1}{c}{Capacity } \\
    \midrule
1 & 12 & 777 & 17 & 6 & 235 & 32 & 10 & 300 & 47 & 10 & 151 & 62 & 12 & 152 \\
2 & 8 & 702 & 18 & 4 & 702 & 33 & 12 & 147 & 48 & 8 & 146 & 63 & 10 & 152 \\
3 & 12 & 777 & 19 & 4 & 147 & 34 & 8 & 146 & 49 & 4 & 157 & 64 & 12 & 152 \\
4 & 10 & 149 & 20 & 6 & 235 & 35 & 8 & 702 & 50 & 6 & 590 & 65 & 4 & 157 \\
5 & 8 & 702 & 21 & 20 & 151 & 36 & 12 & 147 & 51 & 16 & 150 & 66 & 6 & 146 \\
6 & 8 & 513 & 22 & 10 & 151 & 37 & 6 & 777 & 52 & 4 & 157 & 67 & 6 & 288 \\
7 & 8 & 702 & 23 & 10 & 300 & 38 & 6 & 777 & 53 & 4 & 145 & 68 & 10 & 152 \\
8 & 8 & 513 & 24 & 20 & 151 & 39 & 8 & 153 & 54 & 4 & 702 & 69 & 4 & 157 \\
9 & 4 & 533 & 25 & 6 & 417 & 40 & 8 & 146 & 55 & 6 & 590 & 70 & 8 & 150 \\
10 & 12 & 147 & 26 & 6 & 417 & 41 & 10 & 154 & 56 & 8 & 702 & 71 & 8 & 148 \\
11 & 4 & 533 & 27 & 10 & 300 & 42 & 8 & 148 & 57 & 6 & 437 & 72 & 8 & 150 \\
12 & 8 & 148 & 28 & 12 & 405 & 43 & 12 & 405 & 58 & 4 & 145 & 73 & 4 & 152 \\
13 & 10 & 300 & 29 & 8 & 146 & 44 & 10 & 154 & 59 & 8 & 150 & 74 & 8 & 153 \\
14 & 10 & 149 & 30 & 16 & 150 & 45 & 6 & 437 & 60 & 8 & 702 & 75 & 6 & 146 \\
15 & 8 & 148 & 31 & 12 & 147 & 46 & 6 & 288 & 61 & 8 & 150 & 76 & 4 & 152 \\
16 & 4 & 147 \\
    \bottomrule
    \end{tabular}%
    }
  \label{tab:capa_fft_base}%
\end{table}%


\end{document}